\documentclass{amsart}
\usepackage[dvips]{graphicx}
\usepackage{amstext,amssymb}
\numberwithin{equation}{section}
\newtheorem{theorem}{Theorem}[section]

\newtheorem{definition}{Definition }[section]
\newtheorem{lemma}{Lemma}[section]

\newtheorem{remark}{Remark}[section]

\newcommand{\norm}[1]{\left\Vert#1\right\Vert}
\newcommand{\norml}[2]{\left\Vert#1\right\Vert_{L^2(#2)}}

\newcommand{\normeone}[1]{\|{\hskip -0.6pt} | #1 \|{\hskip -0.6pt} |_{1,h,1}}
\newcommand{\normeq}[1]{\|{\hskip -0.6pt} | #1 \|{\hskip -0.6pt} |_{1,h,q}}
\newcommand{\abs}[1]{\bigl\vert#1\bigr\vert}
\newcommand{\pd}[1]{\left\langle #1\right\rangle}
\newcommand{\pdaj}[2]{\left\langle \left\{#1\right\}, \left[#2\right]\right\rangle_e}
\newcommand{\pdja}[2]{\left\langle \left[#1\right], \left\{#2\right\}\right\rangle_e}
\newcommand{\pdjj}[2]{\left\langle \left[#1\right], \left[#2\right]\right\rangle_e}
\newcommand{\set}[1]{\left\{#1\right\}}
\newcommand{\av}[1]{\left\{#1\right\}}
\newcommand{\jm}[1]{\left[#1\right]}

\newcommand{\tuh}{\tilde{u}_h}

\newcommand{\db}{\displaybreak[0]}
\newcommand{\nn}{\nonumber}

\newcommand{\al}{\alpha}
\newcommand{\be}{\beta}

\newcommand{\De}{\Delta}
\newcommand{\ep}{\varepsilon}
\newcommand{\ga}{\gamma}
\newcommand{\Ga}{\Gamma}

\newcommand{\La}{\Lambda}
\newcommand{\na}{\nabla}

\newcommand{\Om}{\Omega}
\newcommand{\pa}{\partial}
\newcommand{\pr}{\prime}
\newcommand{\si}{\sigma}
\newcommand{\ta}{\theta}
\newcommand{\vp}{\varphi}

\newcommand{\csta}[1]{C_{{\rm sta},#1}}

\newcommand{\cerr}[1]{C_{{\rm err},#1}}
\newcommand{\cerrh}[1]{\widehat{C}_{{\rm err},#1}}
\renewcommand{\i}{{\rm\mathbf i}}

\DeclareMathOperator{\re}{{Re}}
\DeclareMathOperator{\im}{{Im}}

\newcommand{\br}{\mathbf{R}}
\newcommand{\bR}{\mathbf{R}}

\newcommand{\p}{\partial}

\newcommand{\Ome}{\Omega}
\newcommand{\nab}{\nabla}
\newcommand{\Del}{\Delta}

\newcommand{\cT}{\mathcal{T}}

\newcommand{\cE}{\mathcal{E}}

\newcommand{\loc}{{\rm loc}}

\newcommand{\Langle}{\left\langle}
\newcommand{\Rangle}{\right\rangle}

\def\jump#1{[#1]}
\def\Jump#1{\left[#1\right]}
\def\avrg#1{\{#1\}}

\begin{document}
\title[$HP$-Discontinuous Galerkin Methods for the Helmholtz Equation]
{$hp$-discontinuous Galerkin Methods for the Helmholtz Equation with Large Wave Number}

\author{Xiaobing Feng}
\address{Department of Mathematics \\
         The University of Tennessee \\
         Knoxville, TN 37996, U.S.A.}
\email{xfeng@math.utk.edu}

\author{Haijun Wu}
\address{Department of Mathematics \\
Nanjing University\\
Nanjing, Jiangsu 210093, P. R. China.}
\email{hjw@nju.edu.cn}

\thanks{The work of the first author was partially supported by the NSF
grants DMS-0410266 and DMS-0710831.  The work of the second author was
partially supported by the national basic research program of China
under grant 2005CB321701 and by the program for the new century outstanding
talents in universities of China. Part of this joint work was carried out
during the second author's recent visit of the University of Tennessee,
the author would like to thank the host institution for its hospitality and
financial support of the visit.}

\keywords{
Helmholtz equation, time harmonic waves, absorbing boundary conditions,
$hp$-discontinuous Galerkin methods, error estimates
}

\subjclass{
65N12, %Stability and convergence of numerical methods
65N15, %Error bounds
65N30, %Finite elements, Rayleigh-Ritz and Galerkin methods, finite methods
78A40  %Wave and radiation
}

\begin{abstract}
This paper develops some interior penalty $hp$-discontinuous
Galerkin ($hp$-DG) methods for the Helmholtz equation in two and three
dimensions. The proposed $hp$-DG methods are defined using a sesquilinear 
form which is not only mesh-dependent but also degree-dependent.
In addition, the sesquilinear form contains penalty terms which not only
penalize the jumps of the function values across the element edges but also
the jumps of the first order tangential derivatives as well as jumps of
all normal derivatives up to order $p$. Furthermore, to ensure
the stability, the penalty parameters are taken as complex numbers
with positive imaginary parts. It is proved that the proposed $hp$-discontinuous 
Galerkin methods are absolutely stable (hence, well-posed).
For each fixed wave number $k$, sub-optimal order error estimates in the 
broken $H^1$-norm and the $L^2$-norm are derived without any mesh constraint.
The error estimates and the stability estimates are improved to optimal 
order under the mesh condition  $k^3h^2p^{-1}\le C_0$ by utilizing these 
stability and error estimates and using a stability-error iterative procedure
To overcome the difficulty caused by strong indefiniteness 
of the Helmholtz problems in the stability analysis
for numerical solutions, our main ideas for stability analysis are to make
use of a local version of the Rellich identity (for the Laplacian)
and to mimic the stability analysis for the PDE solutions given in
\cite{cummings00,Cummings_Feng06,hetmaniuk07}, which enable us to derive
stability estimates and error bounds with explicit
dependence on the mesh size $h$, the polynomial degree $p$,
the wave number $k$, as well as all the penalty parameters
for the numerical solutions.
\end{abstract}

\maketitle

%\setcounter{page}{1}

%%%%%%%%%%%%%%
\section{Introduction}\label{sec-1}
This is the second installment in a series (cf. \cite{fw08a}) which devotes
to developing and analyzing novel interior penalty discontinuous Galerkin
(IPDG) methods for the following Helmholtz problem with large wave number:
\begin{alignat}{2}
-\Del u - k^2 u &=f  &&\qquad\mbox{in  } \Ome:=\Ome_1\setminus D,\label{e1.1}\\
\frac{\pa u}{\pa n_\Ome} +\i k u &=g &&\qquad\mbox{on } \Gamma_R:=\p\Ome_1,
\label{e1.2} \\
u &=0 &&\qquad\mbox{on } \Ga_D:=\p D,\label{e1.3}
\end{alignat}
where $k\in\br$, called wave number, is a (large) positive number,
$D\subset\Ome_1\subset \br^d,\, d=2,3$, $D$ is known as
a scatterer and is assumed to be  a bounded Lipschitz domain,
$\Ome_1$, which is assumed to be a polygonal/polyhedral
domain and often taken as a $d$-rectangle in applications, defines
the size of the computational domain. Note that
$\pa\Ome=\Gamma_R\cup \Gamma_D$. $n_\Ome$ denotes the unit outward
normal to $\Ome$. $\i:=\sqrt{-1}$ denotes the imaginary unit.
Condition \eqref{e1.2} with $g=0$ is known as the
first order absorbing boundary condition (cf. \cite{em79}), which
is used to minimize the reflection at the boundary $\Gamma_R$
and to limit the computation of the original
scattering problem just on the finite domain $\Ome$.
Boundary condition \eqref{e1.3} implies that the scatterer is sound-soft.
We note that the case $D=\emptyset$ also arises in applications
either as a consequence of frequency domain treatment of waves or
when time-harmonic solutions of the scalar wave equation are sought
(cf. \cite{dss94}).

In \cite{fw08a} we proposed and analyzed some IPDG methods for problem
\eqref{e1.1}--\eqref{e1.3} using piecewise linear
polynomial trial and test functions. It was proved that the proposed
methods are unconditionally (with respect to mesh size $h$)
stable and well-posed for all wave number $k>0$. Optimal order
error estimates were established showing explicit
dependence of the error bounds on $h$, $k$ and all penalty parameters.
However, due to the existence of a pollution term, the (broken)
$H^1$-norm error bound deteriorates as the wave number $k$ increases
under the practical ``rule of thumb" mesh constraint $kh\lesssim 1$.
To improve the accuracy and efficiency of those IPDG methods,
it is necessary to use (piecewise) high order polynomial
trial and test functions partly because of the rigidity and low
approximability of linear functions and partly because of the very
oscillatory nature of high frequency waves.
However, simply replacing the linear element by high order
elements in the IPDG methods of \cite{fw08a} does not
reduce the pollution very much, in particular, the
theoretical error bounds do not change much because the analysis
of \cite{fw08a} indeed strongly depends on the properties of
linear functions.

Motivated by the above challenge and observation, the primary
goal of this paper is to develop some new $hp$-interior
penalty discontinuous Galerkin ($hp$-IPDG) methods which
retain the advantages of the IPDG methods of \cite{fw08a} but
improve their accuracy and stability by exploiting the efficiency
and flexibility of piecewise high order polynomial functions.
To the end, our key idea is to construct a sesquilinear form (as a
discretization of the Laplacian) which is not only
mesh-dependent (or $h$-dependent) but also degree-dependent
(or $p$-dependent) by introducing penalty terms which not only
penalize the jumps of the function values across the element edges but also
the jumps of the first order tangential derivatives as well as jumps of
all normal derivatives up to order $p$. In addition, as in \cite{fw08a},
to ensure the stability, all penalty parameters are taken as complex
numbers with positive imaginary parts.
Since the Helmholtz equation with large wave number is
non-Hermitian and strongly indefinite, as expected, stability
estimates (or a priori estimates) for numerical solutions
under practical mesh constraints is a difficult task
to carry out regardless which discretization method is used.
To overcome the difficulty, as in \cite{fw08a}, the cruxes of
our analysis are to establish and to make use of a local
version of the Rellich identity (for the Laplacian) and to mimic
the stability analysis for the PDE solutions given in
\cite{cummings00,Cummings_Feng06,hetmaniuk07}. The key idea here
is to use the special test function $\nab u_h\cdot (x-x_\Ome)$
(defined element-wise) with $u_h$ denoting the $hp$-IPDG solution, such a
test function is valid for any DG method. We remark that the same
technique was successfully employed by Shen and Wang in \cite{sw07} to
establish the stability and error analysis for the spectral Galerkin
approximation of the Helmholtz problem. We also note that
although the similar techniques to those in \cite{fw08a,sw07}
are utilized in this paper to carry out the stability analysis,
the analysis of this paper is more involved
because the special sesquilinear form of this paper, which contains jumps
of high order normal derivatives, is a lot more complicate to deal with,
even they are similar conceptually.

Since the Helmholtz equation appears, in one way or another, directly
or indirectly, in almost all wave-related problems arisen from many science,
engineering, and industry applications, solving the Helmholtz equation,
in one form or another, has always been and remains at the center of
wave computation. We refer the reader to
(\cite{amm06,ak79,aw80,bao95,chang90,ce06,cm87b,cummings00,
dss94,er03,goldstein81,hh92,ib97,
perugia07,ked08,zienkiewicz00} and the references therein)
for some recent developments on numerical methods, in particular,
Galerkin type methods, for the Helmholtz equation. We also
refer the reader to \cite{fw08a} for a brief review about
some theoretical issues for finite element approximations
(and other types of Galerkin approximations) of the Helmholtz equation.

The $hp$-finite element method ($hp$-FEM) is a modern version of the
finite element method, capable of achieving exceptionally fast (exponential)
convergence. It combines the flexibility of the standard finite element method
and the high order accuracy of the spectral method. Consequently, the $hp$-FEM
can often attain more accurate results than the standard finite element
method does while using less CPU time and resources. The $hp$-FEM has
undergone intensive developments both on theory and implementation
in the past twenty five years. We refer the reader to the survey paper
\cite{bg96} and two recent monographs \cite{schwab98,ssd04} for
a detailed exposition on the basic theory and advanced topics of the $hp$-FEM.

Discontinuous Galerkin (DG) methods was first proposed in 1970s,
they were not popular then because they produce larger
algebraic systems than standard finite element methods do.
However, due to the emergence of high performance computers and fast
solvers since early 1990s, especially, massively parallel computers
and parallel solvers such as multilevel and domain decomposition
methods, which together with advantages of DG methods has
quickly attracted renewed interests in DG methods. They have been
heavily developed and tested in the past fifteen years, we refer the
reader to \cite{abcm01} and the references therein for a review of
recent developments. As is well known now,
DG methods have several advantages over other types of numerical methods.
For example, the trial and test spaces are easy to construct, they can
naturally handle inhomogeneous boundary conditions and curved boundaries;
they also allow the use of highly nonuniform and unstructured meshes,
and have built-in parallelism which permits coarse-grain parallelism.
In addition, the fact that the mass matrices are block diagonal is an
attractive feature in the context of time-dependent problems, especially
if explicit time discretizations are used. Moreover, as proved in
\cite{fw08a}, DG methods are also effective and have advantages over
finite element methods for the strongly indefinite Helmholtz equation,
which has not been well understood before.  We refer the
reader to \cite{arnold82,abcm01,baker77,CKS00,CS98,dd76,fk07,rwg99,w78}
and the references therein for a detailed account on DG methods for
coercive elliptic and parabolic problems, and to \cite{btp07,ews06,
gs06,hjs02,hss02} and the references therein for recent developments
on $hp$-discontinuous Galerkin ($hp$-DG) methods.

The remainder of this paper is organized as follows. In Section \ref{sec-2},
we first introduce notation and gather some preliminaries,
and then formulate our $hp$-IPDG methods.
Both symmetric and non-symmetric methods are constructed
and various possible variants are also discussed.
Section \ref{sec-sta} devotes to the stability analysis for the $hp$-IPDG
methods proposed in Section \ref{sec-2}. It is proved that the proposed
$hp$-IPDG methods are stable (hence well-posed) without any mesh constraint.
In Section \ref{sec-err}, using the stability results of
Section \ref{sec-sta} we prove that for each fixed wave number $k$,
sub-optimal order (with respect to $h$ and $p$)
error estimates in the broken $H^1$-norm and the $L^2$-norm are derived without any mesh constraint.
Finally, using the stability estimate of Section \ref{sec-sta},
the error estimates of Section \ref{sec-err} and a
stability-error iterative procedure we obtain some much improved (optimal order) stability
and error estimates for the $hp$-IPDG solutions under the mesh condition
 $k^3h^2p^{-1}\le C_0$ in Section \ref{sec-SEII}, where $C_0$ is some constant independent of $k$, $h$, $p$, and the penalty parameters.

%%%%%%%%%%%%%%
\section{Formulation of $hp$-interior penalty discontinuous Galerkin methods}
\label{sec-2}

\subsection{Notation and preliminaries}
The space, norm and inner product notation used in this paper all are
standard, we refer to \cite{bs94,ciarlet78,baker77} for their precise
definitions. On the other hand, we note that all functions in this paper
are complex-valued, so the familiar terminologies such ``symmetric/non-symmetric"
and ``bilinear" are replaced respectively by terms ``Hermitian/non-Hermitian"
and ``sesquilinear".  For a complex number $a=a_r+\i a_i$ ($a_r$ and $a_i$
are real numbers), $\overline{a}:= a_r -\i a_i$ denotes the complex conjugate of $a$.
$(\cdot,\cdot)_Q$ and $\langle \cdot,\cdot\rangle_\Sigma$
for $\Sigma\subset \pa Q$ denote the complex $L^2$-inner product
on $L^2(Q)$ and $L^2(\Sigma)$
spaces, respectively. $(\cdot,\cdot):=(\cdot,\cdot)_\Ome$
and $\langle \cdot,\cdot\rangle:=\langle \cdot,\cdot\rangle_{\p\Ome}$.
We also define
\[
H_{\Ga_D}^1(\Om):=\set{u\in H^1(\Om);\, u=0 \text{ on }\Ga_D}.
\]
Throughout the paper, $C$ is used to denote a generic positive constant
which is independent of  $k$, $h$, $p$, and the penalty parameters. We also use the shorthand notation
$A\lesssim B$ and $B\gtrsim A$ for the inequality $A\leq C B$ and $B\geq CA$.
$A\simeq B$ is for the statement $A\lesssim B$ and $B\lesssim A$.

We now give the definition of star-shaped domains.

\begin{definition}\label{def1}
$Q\subset \bR^d$ is said to be a {\em star-shaped} domain with respect
to $x_Q\in Q$ if there exists a nonnegative constant $c_Q$ such that
\begin{equation}\label{estar}
(x-x_Q)\cdot n_Q\ge c_Q \qquad \forall x\in\pa Q.
\end{equation}
$Q\subset \bR^d$ is said to be {\em strictly star-shaped} if $c_Q$ is positive.
\end{definition}

In this paper, we assume that $\Ome_1$ is a strictly star-shaped domain.
Recall that $\Ome_1$ is often taken as a $d$-rectangle in practice.
We also assume that the scatterer $D$ is a star-shaped domain, without loss
of the generality, with respect to the same point $x_{\Ome_1}$ as $\Ome_1$
does.  This implies that $x_{\Ome_1}\in D\subset \Ome_1$.
Under these assumptions, there hold following stability estimates
for problem \eqref{e1.1}--\eqref{e1.3}.
\begin{theorem}\label{stability}
Suppose $\Om_1\subset \bR^d$ is a strictly star-shaped domain
and $D\subset\Om_1$ is a star-shaped domain. Then the solution
$u$ to problem \eqref{e1.1}--\eqref{e1.3} satisfies
\begin{eqnarray}\label{e2.1}
\|u\|_{H^j(\Ome)} \lesssim \Bigl(\frac{1}k+ k^{j-1} \Bigr)
\bigl( \|f\|_{L^2(\Ome)} + \|g\|_{L^2(\Ga_R)} \bigr)
\end{eqnarray}
for $j=0, 1$ if $u\in H^{3/2+\ep}(\Om)$ for some $\ep>0$. \eqref{e2.1}
also holds for $j=2$ if $u\in H^2(\Om)$. Furthermore, there hold
\begin{align}
\|u\|_{H^j_{\mbox{\small loc}}(\Ome)} &\lesssim \Bigl(\frac{1}k+ k^{j-1} \Bigr)
\Bigl( \|f\|_{L^2(\Ome)} + \|g\|_{L^2(\Gamma_R)}
+\sum_{\ell=1}^{j-2} \norm{f}_{H^\ell_{\small loc}(\Ome)}  \Bigr) \label{e2.1a}
\end{align}
for $j=3,4,\cdots,q$ if $u\in H^2(\Ome)\cap H^q_{\mbox{\small loc}}(\Ome)$
for some positive integer $q\geq 3$.
\end{theorem}

\begin{proof}
Inequality \eqref{e2.1} for $j=0,1,2$ was proved in
\cite{cummings00,Cummings_Feng06,hetmaniuk07}.
Inequality \eqref{e2.1a} follows from \eqref{e2.1} and an application
of the standard cutoff function technique together
with an induction argument. We leave the derivation to the interested reader.
\end{proof}

%%%%%%%%%%%%%%
\subsection{Formulation of $hp$-IPDG methods}
To formulate our $hp$-IPDG methods, we need to introduce some notation,
most of them were already appeared in \cite{fw08a}.
Let $\cT_h$ be a family of partitions of the domain $\Ome:=\Ome_1\setminus D$
parameterized by $h\in (0,h_0)$. For any $K\in \cT_h$, we define
$h_K:=\mbox{diam}(K)$.
Similarly, for each edge/face $e$ of $K\in \cT_h$, $h_e:=\mbox{diam}(e)$.
We impose the following mild restrictions on the partition $\cT_h$:
\begin{itemize}
\item[(i)] The elements of $\cT_h$ satisfy the minimal angle condition,
\item[(ii)] $\cT_h$ is locally quasi-uniform, that is if two elements $K$
and $K'$ are adjacent ( i.e., $\mbox{meas}(\p K\cap\p K') > 0$ ), then
$h_{K} \simeq h_{K'}$. Where $\mbox{meas}(e)$ stands for $(d-1)$-dimensional
Lebesgue measure of $e$.
\end{itemize}
For convenience, we assume $\mathrm{diam} (\Om)\simeq 1$,
hence $h_e, h_K\lesssim 1.$

For any two elements $K,\, K'\in\cT_h$, we call
$e=\p K\cap \p K'$ an {\em interior edge/face} of
$\cT_h$ if $\mbox{meas}(e) >0$. Note that $e$ could be
portion of a side/face of the element $K$ or $K'$ in the case of
geometrically nonconforming partition. Also, for any element $K\in\cT_h$,
we call $e=\p K\cap \p\Ome$ a {\em boundary edge/face} if
$\mbox{meas}(e) >0$. Then we define
\begin{align*}
\cE_h^I&:= \mbox{ set of all interior edges/faces of } \cT_h,\\
\cE_h^R&:= \mbox{ set of all boundary edges/faces of $\cT_h$ on $\Ga_R$} ,\\
\cE_h^D&:= \mbox{ set of all boundary edges/faces of $\cT_h$ on $\Ga_D$} ,\\
\cE_h^{RD}&:= \cE_h^R\cup\cE_h^D= \mbox{ set of all boundary edges/faces of $\cT_h$} ,\\
\cE_h^{ID}&:=\cE_h^I\cup\cE_h^D= \mbox{ set of all edges/faces of $\cT_h$ except those on $\Ga_R$},\\
\cE_h&:=\cE_h^I\cup \cE_h^{RD}=\mbox{ set of all edges/faces of } \cT_h.
\end{align*}
We also define the jump $\jump{v}$ of $v$ on an interior edge/face
$e=\p K\cap \p K'$ as
\[
\jump{v}|_{e}:=\left\{\begin{array}{ll}
       v|_{K}-v|_{K'}, &\quad\mbox{if the global label of $K$ is larger},\\
       v|_{K'}-v|_{K}, &\quad\mbox{if the global label of $K'$ is larger}.
\end{array} \right.
\]
If $e\in\cE_h^D$, set $\jump{v}|_{e}=v|_{e}$. The following convention is
adopted in this paper
\[
\avrg{v}|_{e}:=\frac12\bigl( v|_{K}+ v|_{K'} \bigr)\qquad \mbox{if }
e=\p K\cap \p K'.
\]
If $e\in\cE_h^{RD}$, set $\avrg{v}|_{e}=v|_{e}$.  
For every $e=\p K\cap \p K^\pr\in\cE_h^I$, let $n_e$ be the unit outward normal
to edge/face $e$ of the element $K$ if the global label of $K$ is bigger
and of the element $K^\pr$ if the other way around. For every $e\in\cE_h^{RD}$, let $n_e=n_\Om$ the unit outward normal to $\pa\Om$.

Let $p\geq 1$ be a fixed integer, which will be used to
denote the degree of the $hp$-IPDG methods in this paper.
For each integer $0\le q\le p$, we define the ``energy" space
\[
E^q:=\prod_{K\in\cT_h} H^{q+1} (K),
\]
and the sesquilinear form $a^q_h(\cdot,\cdot)$ on $E^q\times E^q$
\begin{equation}\label{eah}
a_h^q(u,v) :=b_h(u,v) + \i\Bigl( L_1(u,v)+ \sum_{j=0}^q J_j(u,v) \Bigr) \qquad
\forall u, v\in E^q,
\end{equation}
\begin{align} \label{ebh}
b_h(u,v) &:=\sum_{K\in\cT_h} (\nab u,\nab v)_K
-\sum_{e\in\cE_h^{ID}} \left( \pdaj{\frac{\pa u}{\pa n_e}}{v} \right. \\
&\hskip 1.2in \left. +\si\pdja{u}{\frac{\pa v}{\pa n_e}}\right),\nonumber \\
L_1(u,v) &:=\sum_{e\in\cE_h^{ID}}\sum_{\ell=1}^{d-1}\frac{\be_{1,e} p}{h_e} \pdjj{\frac{\pa u}{\pa \tau_e^\ell}}
{\frac{\pa v}{\pa \tau_e^\ell}},\label{eL1}\\
J_0(u,v)&:=\sum_{e\in\cE_h^{ID}}\frac{\ga_{0,e}\,p}{h_e} \pdjj{u}{v},\label{eJ0}\db\\
J_j(u,v) &:=\sum_{e\in\cE_h^I} \ga_{j,e} \left(\frac{h_e}{p}\right)^{2j-1}\pdjj{\frac{\p^j u}{\p n_e^j}}{\frac{\p^j v}{\p n_e^j}},
\qquad j=1,2,\cdots, q,
\label{eJ1}
\end{align}
and $\si$ is a real number. $\ga_{0,e},\cdots, \ga_{q,e}>0$ and $\be_{1,e}\ge 0$
are  numbers to be specified later. $\{\tau^\ell_e\}_{\ell=1}^{d-1}$
denote an orthogonal coordinate frame on the edge/face $e\in \cE_h$,
$\frac{\pa u}{\pa \tau_e^\ell} :=\nab u\cdot \tau_e^\ell$ stands for the
tangential derivative of $u$ in the direction $\tau_e^\ell$, and
$\frac{\p^j u}{\p n_e^j}$ denotes the $j$th order
normal derivative of $u$ on $e$.

It is easy to check that $(-\Del u,v)=a_h^q(u,v)$ for all $u\in H^{q+1}(\Ome)$
and $v\in E^q$. Hence, $a_h^q(\cdot,\cdot)$ is a consistent discretization for
$-\Del$. When $\si=1$,  $a_h^q(\cdot,\cdot)$ is symmetric, that is,
$a_h^q(u,v)=a_h^q(\overline{v},\overline{u})$.
On the other hand, when $\si\neq 1$,
$a_h^q(\cdot,\cdot)$ is non-symmetric. In particular, $\si =-1$
would correspond to the non-symmetric IPDG method studied in \cite{rwg99}
for coercive elliptic problems.
In this paper, for the ease of presentation, we only consider the case $\si=1$.
The penalty constants in $\i\bigl(L_1(u,v)+J_0(u,v)+\cdots+J_q(u,v)\bigr)$ are
$\i\beta_{1,e},\i\gamma_{0,e},\cdots,\i\gamma_{q,e}$, respectively.  So they are
pure imaginary numbers with positive imaginary parts. It turns out that if any
of them is replaced by a complex number with positive
imaginary part, the ideas of the paper still apply.  Here we set their
real parts to be zero because the terms from real parts do not help much
(and do not cause any problem either) in our analysis.

Next, we introduce the following semi-norms on the space $E^q$:
\begin{align}
\abs{v}_{1,h}&:=\Bigl(\sum_{K\in\cT_h} \norml{\nab v}{K}^2\Bigr)^{\frac12},
\label{e2.5b}\db\\
\norm{v}_{1,h,q}&:=\left(\abs{v}_{1,h}^2
+\sum_{e\in\cE_h^{ID}} \left(\frac{\ga_{0,e}\,p}{h_e}\norml{\jm{v}}{e}^2
+\sum_{\ell=1}^{d-1}\frac{\be_{1,e}p}{h_e}\norml{\jm{\frac{\pa v}{\pa \tau_e^\ell}}}{e}^2\right) \right. \label{e2.5} \\
&\hskip 0.75in \left.
+ \sum_{j=1}^q\sum_{e\in\cE_h^I} \ga_{j,e} \left(\frac{h_e}{p}\right)^{2j-1}
\norml{\jm{\frac{\pa^j v}{\pa n_e^j}}}{e}^2 \right)^{\frac12},  \nonumber\\
\normeq{v}&:=\left(\norm{v}_{1,h,q}^2+\sum_{e\in\cE_h^{ID}}
\frac{h_e}{\ga_{0,e}\,p}\norml{\av{\frac{\pa v}{\pa n_e}}}{e}^2 \right)^{\frac12}.
\label{e2.5a}
\end{align}
Clearly, $\norm{\cdot}_{1,h,q}$ and $\normeq{\cdot}$ are norms on $E^q$
if $\p D\neq\emptyset$ but only semi-norms if $\p D=\emptyset$.

It is easy to check that the sesquilinear form $a_h^q(\cdot,\cdot)$
satisfies: For any $v\in E^q$
\begin{align}
\re a_h^q(v,v)
&=\abs{v}_{1,h}^2-2\re \sum_{e\in\cE_h^{ID}}\pdaj{\frac{\pa v}{\pa n_e}}{v},
\label{ea1}\\
\im a_h^q(v,v)&=L_1(v,v)+J_0(v,v)+\cdots+J_q(v,v). \label{ea2}
\end{align}

Using the sesquilinear form $a_h^q(\cdot,\cdot)$ we now
introduce the following weak formulation for \eqref{e1.1}--\eqref{e1.2}:
Find $u\in E^q\cap H^1_{\Ga_D}(\Ome)\cap H^2_{\loc}(\Ome)$ such that
\begin{equation}
a_h^q(u,v) - k^2(u,v) +\i k \langle u,v\rangle_{\Gamma_R}
=(f,v)+\langle g, v \rangle_{\Gamma_R}
\label{2.6}
\end{equation}
for any $v\in E^q\cap H^1_{\Ga_D}(\Ome)\cap H^2_{\loc}(\Ome)$.
The above formulation is consistent with the boundary value problem
\eqref{e1.1}--\eqref{e1.2} because $a_h^q(\cdot,\cdot)$ is consistent
with $-\Del$.

For any $K\in \cT_h$, let $\mathcal{P}_p(K)$ denote the set of all polynomials
whose degrees do not exceed $p$. We define our $hp$-IPDG approximation
space $V_h^p$ as
\[
V_h^p:=\prod_{K\in \cT_h} \mathcal{P}_p(K).
\]
Clearly, $V_h^p\subset E^q\subset L^2(\Ome)$.
But $V_h^p\not\subset H^1(\Ome)$.  We are now ready to define
our $hp$-IPDG methods based on the weak formulation
\eqref{2.6}:
For each $0\leq q\leq p$, find $u_h^q\in V_h^p$ such that
\begin{equation}\label{edg}
a_h^q(u_h^q, v_h) - k^2(u_h^q, v_h)
+\i k \langle u_h^q, v_h\rangle_{\Gamma_R}
=(f,  v_h)+\langle g,  v_h \rangle_{\Gamma_R}  \qquad\forall  v_h\in V_h^p.
\end{equation}

\begin{remark}
(a) When $p=q=1$, the above method \eqref{edg} is exactly the scheme
proposed in \cite{fw08a}. The $L_1$ term, which
penalizes the jumps of the first order tangential derivatives,
plays an important role for getting a better (theoretical) stability
estimate in \cite{fw08a}. However, our analysis to be given in the next section
suggests that the $L_1$ term plays a less pivotal role
for high order IPDG methods.

(b) In fact, \eqref{edg} defines $p+1$ different IPDG methods for $q=0,1,\cdots,p$.
$q=1$ would correspond to using high order elements in the IPDG formulation
proposed in \cite{fw08a}.

(c) The idea of penalizing the jumps of normal derivatives (i.e., the $J_1$
term above) for second  order PDEs was used early by Douglas and
Dupont \cite{dd76}
in the context of $C^0$ finite element methods, by Baker \cite{baker77}
(with a different weighting, also see \cite{fk07}) for fourth order PDEs.
The idea of using multipenalties $J_0, J_1, \cdots, J_p$ with positive
penalty parameters was first used by Arnold in \cite{arnold82} for coercive
elliptic and parabolic PDEs. The use of $L_1$ term was first introduced
in \cite{fw08a}.
\end{remark}

In the next two sections, we shall study the stability and error analysis
for the $hp$-IPDG method \eqref{edg}. 
Especially, we are interested in knowing how the stability constants and
error constants depend on the wave number $k$ (and mesh size $h$ and
element degree $p$, of course) and on the penalty parameters,
and what are the ``optimal" relationship
between mesh size $h$ and the wave number $k$.

%%%%%%%%%%%%%%%%%%%%%%%%%%%%%%%%%%%%%%%%%%
\section{Stability estimates}\label{sec-sta}

The goal of this section is to derive stability estimates (or a priori estimates)
for schemes \eqref{edg}. To the end, momentarily, we assume that the solution $u_h^q$
to \eqref{edg} exists and will revisit the existence and uniqueness
issues later at the end of the section. We like to note that because its strong
indefiniteness, unlike in the case of coercive elliptic and parabolic problems
(cf. \cite{arnold82,abcm01,baker77,dd76,fk07,rwg99,w78}),
the well-posedness of scheme \eqref{edg} is difficult to prove
under practical mesh constraints.

To derive stability estimates for scheme \eqref{edg}, our approach is to mimic
the stability analysis for the Helmholtz problem \eqref{e1.1}--\eqref{e1.2}
given in \cite{cummings00,Cummings_Feng06,hetmaniuk07}.
The key ingredients of our analysis
are to use a special test function $v_h=\al\cdot\na u_h^q$ (defined element-wise)
with $\al(x):=x-x_{\Ome_1}$ in \eqref{edg} and to use the Rellich
identity (cf. \cite{Cummings_Feng06} and below) on each element.
Due to existence of multiple penalty terms in $a_h^p(\cdot,\cdot)$,
which do not appear in \cite{cummings00,Cummings_Feng06,hetmaniuk07}, the analysis
to be given below is much more delicate and complicate than those of
\cite{cummings00,Cummings_Feng06,hetmaniuk07},
although they are similar conceptually.  Since most proofs of this
section are in the same lines as those of the proofs
in Section 4 of \cite{fw08a}, we shall omit some details if they are
already given in \cite{fw08a}, but shall provide them if there
are meaningful differences.

We first cite the following lemma which establishes three integral identities
and play a crucial role in our analysis.  A proof of the lemma
can be found in \cite[Lemma 4.1]{fw08a}.

\begin{lemma}\label{lem2.1}
Let $\al(x):=x-x_{\Ome_1}$, $v\in E^1$, $K, K'\in\cT_h$ and $e\in \cE_h^{ID}$.
Then there hold
\begin{align}
&d\norml{v}{K}^2+2\re(v,\al\cdot\na v)_K
=\int_{\pa K}\al\cdot n_K\abs{v}^2,\label{eid1}\\
&(d-2)\norml{\na v}{K}^2+2\re\big(\na v,\na(\al\cdot\na v)\big)_K
=\int_{\pa K}\al\cdot n_K\abs{\na v}^2,\label{eid2}\\
&\pdaj{\frac{\pa v}{\pa n_e}}{\al\cdot\na v}
-\pd{\al\cdot n_e\av{\na v},\jm{\na v}}_e  \label{eid3}\\
&\hskip 1.35in=\sum_{\ell=1}^{d-1}\int_e\left(\al\cdot\tau_e^\ell\av{\frac{\pa v}{\pa n_e}}
-\al\cdot n_e\av{\frac{\pa v}{\pa \tau_e^\ell}}\right)\jm{\frac{\pa\overline{v}}{\pa\tau_e^\ell}}, \nn
\end{align}
where $x_{\Ome_1}$ denotes the point in the star-shaped domain definition for
$\Ome_1$ (see Definition \ref{def1}).
\end{lemma}

\begin{remark}
The identity \eqref{eid2} can be viewed as a local version of the Rellich
identity for the Laplacian $\Del$ (cf. \cite{cummings00,Cummings_Feng06}).
Since $V_h^p\subset E^1$, hence, \eqref{eid1}--\eqref{eid3} hold for
any function $v=v_h\in V_h^p$.
\end{remark}

We also need the following trace and inverse inequalities (cf. \cite{schwab98,wh03,be07}).
\begin{lemma}\label{lem-inv}
For any $K\in\cT_h$ and $z\in\mathcal{P}_p(K),$
\begin{align*}
\norml{z}{\pa K}&\lesssim p\,h^{-\frac12}\norml{z}{K}, \\
\norml{\na z}{K}&\lesssim p^2\,h^{-1}\norml{z}{K}.
\end{align*}
\end{lemma}

Now, taking $v_h=u_h^q$ in \eqref{edg} yields
\begin{equation}\label{euh}
a_h^q(u_h^q,u_h^q) -k^2 \norml{u_h^q}{\Om}^2 + \i k \norml{u_h^q}{\Gamma_R}^2
=(f, u_h^q)+\langle g, u_h^q \rangle_{\Gamma_R}.
\end{equation}
Therefore, taking real part and imaginary part of the above equation
and using \eqref{ea1} and \eqref{ea2} we get the following lemma.

\begin{lemma}\label{lem3.1}
Let $u_h^q\in V_h^p$ solve \eqref{edg}. Then
\begin{align}\label{e3.1}
&\abs{u_h^q}_{1,h}^2-2\re \sum_{e\in\cE_h^{ID}}\pdaj{\frac{\pa u_h^q}{\pa n_e}}{u_h^q} -k^2 \norml{u_h^q}{\Om}^2 \\
&\hskip 2.3in
\leq \abs{(f,u_h^q)+\langle g, u_h^q \rangle_{\Gamma_R}}, \nonumber\\
&\sum_{e\in\cE_h^{ID}} \left(\frac{\ga_{0,e} p}{h_e} \norml{\jm{u_h^q}}{e}^2
+\sum_{\ell=1}^{d-1}\frac{\be_{1,e} p}{h_e}
\norml{\jm{\frac{\pa u_h}{\pa \tau_e^\ell}}}{e}^2 \right)
+k\norml{u_h^q}{\Gamma_R}^2 \label{e3.2} \\
&\hskip .4in
+\sum_{j=1}^q \sum_{e\in\cE_h^I} \ga_{j,e} \left(\frac{h_e}{p}\right)^{2j-1}
\norml{\jm{\frac{\pa^j u_h^q}{\pa n_e^j}}}{e}^2
\leq \bigl|(f,u_h^q)+\langle g, u_h^q \rangle_{\Gamma_R}\bigr|.  \nn
\end{align}
\end{lemma}

From \eqref{e3.1} and \eqref{e3.2} we can bound $\abs{u_h^q}_{1,h}$ and
the jumps in terms of $\norml{u_h^q}{\Om}^2$. In order to get the desired a priori
estimates, we need to derive a reverse inequality whose coefficients
can be controlled.  Such a reverse inequality, which is often difficult to get
under practical mesh constraints, and stability estimates for $u_h^q$ will
be derived next.

\begin{theorem}\label{thm_sta}
Let $u_h^q\in V_h^p$ solve \eqref{edg} and suppose
$\beta_{1,e}\geq 0,\ga_{0,e},\cdots,\ga_{q,e}>0$. Then
\begin{align}\label{e3.3}
\norml{u_h^q}{\Om} &+\frac{1}{k} \norm{u_h^q}_{1,h,q} +\norml{u_h^q}{\Gamma_R}
+\frac{1}{k} \Bigl( c_{\Om_1}\sum_{e\in\cE_h^R}\norml{\na u_h^q}{e}^2
\Bigr)^{\frac12} \\
&+\frac{1}{k} \Bigl(\sum_{e\in\cE_h^D} c_D
\Bigl( k^2\norml{u_h^q}{e}^2+\norml{\na u_h^q}{e}^2\Bigr) \Bigr)^{\frac12}
\lesssim \csta{q}\, M(f,g), \nn
\end{align}
where
\begin{align} \label{Mfg}
&M(f,g) :=\norml{f}{\Om}+\|g\|_{L^2(\Gamma_R)},\db\\
&\csta{q}:=\frac1k+\frac{1}{k^2}
+\frac{1}{k^2}\max_{e\in\cE_h^D}\Big(\frac{\ga_{0,e}\,p}{h_e}+\frac{p^5}{\ga_{0,e}h_e^2}+\frac{\beta_{1,e}\,p^5}{h_e^3}+\frac{p^2}{h_e} \Big)\label{csta}\\
&\qquad\qquad+\left\{\begin{aligned}
                  &\frac{1}{k^2}\max_{e\in\cE_h^I}\bigg(\frac{p\,k^2h_e^2+p^5}{\ga_{0,e}\,h_e^2}
+\frac{p}{h_e} \max_{0\le j\le q-1}\sqrt{\frac{\ga_{j,e}}{\ga_{j+1,e}}}
+\frac{p^2}{h_e}&\\
&\hskip 1.1in+\frac{p^3}{h_e^2}\sqrt{\frac{\beta_{1,e}}{\ga_{1,e}}}+\frac{\ga_{q,e}\,p^{2q+3}}{h_e^2} \bigg),\quad\text{ if } q<p, \\
                &\frac{1}{k^2}\max_{e\in\cE_h^I}\bigg(\frac{p\,k^2h_e^2+p^5}{\ga_{0,e}\,h_e^2}
+\frac{p}{h_e} \max_{0\le j\le q-1}\sqrt{\frac{\ga_{j,e}}{\ga_{j+1,e}}}
+\frac{p^2}{h_e}&\\
&\hskip 1.8in+\frac{p^3}{h_e^2}\sqrt{\frac{\beta_{1,e}}{\ga_{1,e}}}\;\bigg),   \quad \text{ if } q=p.
               \end{aligned}\right.
\nn
\end{align}
\end{theorem}

\begin{proof}
We divide the proof into three steps.

{\em Step 1: Derivation of a representation identity for $\|u_h^q\|_{L^2(\Ome)}$.}
Define $v_h\in E$ by $v_h|_K=\al\cdot\na u_h^q|_K$ for every $K\in \cT_h$.
Since $v_h|_K$ is a polynomial of degree no more than $p$ on $K$,
hence, $v_h\in V_h^p$.  Using this $v_h$ as a test function in \eqref{edg}
and taking the real part of the resulted equation we get
\begin{equation}\label{added2}
-k^2 \re (u_h^q, v_h) = \re \bigl( (f,v_h)+\langle g, v_h\rangle_{\Gamma_R}
-a_h^p(u_h^q,v_h)-\i k\pd{u_h^q,v_h}_{\Gamma_R} \bigr) .
\end{equation}

It follows from \eqref{eid1}, \eqref{euh}, and \eqref{added2} that
(compare with (4.11) of \cite{fw08a})
\begin{align}\label{e3.5}
&2k^2\norml{u_h^q}{\Om}^2
=k^2\sum_{K\in\cT_h}\int_{\pa K}\al\cdot n_K\abs{u_h^q}^2
+(d-2)\re \bigl( (f,u_h^q) + \langle g, u_h^q\rangle_{\Gamma_R} \\
&\qquad -a_h^p(u_h^q,u_h^q) \bigr)
+2\re \bigl( (f,v_h)+ \langle g, v_h\rangle_{\Gamma_R}
-a_h^p(u_h^q,v_h)-\i k\pd{u_h^q,v_h}_{\Gamma_R} \bigr)  \nn\db\\
&= k^2\sum_{K\in\cT_h}\int_{\pa K}\al\cdot n_K\abs{u_h^q}^2
+(d-2)\re \bigl( (f,u_h^q) + \langle g, u_h^q\rangle_{\Gamma_R} \bigr) \nn \\
&\qquad +2\re \bigl( (f,v_h) + \langle g, v_h\rangle_{\Gamma_R} \bigr)
+ 2k\im\pd{u_h^q,v_h}_{\Gamma_R} \nn \\
&\qquad -\sum_{K\in\cT_h}\left((d-2)\norml{\na u_h^q}{K}^2
+2\re(\na u_h^q,\na v_h)_K\right)\nn\\
&\qquad +2\sum_{e\in\cE_h^{ID}}\left((d-2)\re\pdaj{\frac{\pa u_h^q}{\pa n_e}}{u_h^q}
+\re\pdaj{\frac{\pa u_h^q}{\pa n_e}}{v_h}\right.\nn\\
&\qquad \left.+\re\pdja{u_h^q}{\frac{\pa v_h}{\pa n_e}} \right)
+2\im\Bigl(L_1(u_h^q,v_h)+\sum_{j=0}^q J_j(u_h^q,v_h) \Bigr). \nn
\end{align}
Using the identity $\abs{a}^2-\abs{b}^2=\re (a+b)(\bar a-\bar b)$ we have
\begin{align}\label{e3.6}
\sum_{K\in\cT_h}\int_{\pa K}\al\cdot n_K\abs{u_h^q}^2
=2\sum_{e\in\cE_h^I}\re \Langle \al\cdot n_e \av{u_h^q},\jm{u_h^q} \Rangle_e
+\Langle \al\cdot n_\Ome, |u_h^q|^2 \Rangle_{\p \Ome}.
\end{align}
Using the identity again followed by the Rellich identity \eqref{eid2} we get
(compare with (4.13) of \cite{fw08a})
\begin{align}\label{e3.7}
&\sum_{K\in\cT_h}\Bigl((d-2)\norml{\na u_h^q}{K}^2
+2\re(\na u_h^q,\na v_h)_K\Bigr)
=\sum_{K\in\cT_h}\int_{\pa K}\al\cdot n_K\abs{\na u_h^q}^2 \\
&\qquad
=2\sum_{e\in\cE_h^I} \re  \Langle \al\cdot n_e\av{\na u_h^q},\jm{\na u_h^q} \Rangle_e
+\sum_{e\in\cE_h^{RD}} \Langle \al\cdot n_e, |\na u_h^q|^2 \Rangle_e \nn \\
&\qquad
=2\sum_{e\in\cE_h^{ID}} \re  \pd{\al\cdot n_e\av{\na u_h^q},\jm{\na u_h^q}}_e
+\sum_{e\in\cE_h^{R}} \pd{\al\cdot n_e, |\na u_h^q|^2}_e \nn \\
&\hskip 2.45in
-\sum_{e\in\cE_h^{D}} \pd{\al\cdot n_e, |\na u_h^q|^2}_e.\nn
\end{align}
Plugging \eqref{e3.6} and \eqref{e3.7} into \eqref{e3.5} gives
(compare with (4.15) of \cite{fw08a})
\begin{align}\label{e3.8}
&2k^2\norml{u_h^q}{\Om}^2\\
&=(d-2)\re\bigl( (f,u_h^q) +\langle g, u_h^q\rangle_{\Gamma_R} \bigr)
+2\re \bigl( (f,v_h)  +\langle g, v_h\rangle_{\Gamma_R} \bigr) \nn\\
&\quad +2k^2\sum_{e\in\cE_h^I}\re\pd{\al\cdot n_e\av{u_h^q},\jm{u_h^q}}_e
+k^2\pd{\al\cdot n_\Ome, |u_h^q|^2}_{\pa\Om}\nn\\
&\quad +2k\im\pd{u_h^q,v_h}_{\Gamma_R}
-\sum_{e\in\cE_h^R}\pd{\al\cdot n_e, |\na u_h^q|^2}_e
+\sum_{e\in\cE_h^D}\pd{\al\cdot n_e, |\na u_h^q|^2}_e
 \nn\\
&\quad -2\sum_{e\in\cE_h^{ID}}\re\pdaj{\frac{\pa u_h^q}{\pa n_e}}{u_h^q}
+2(d-1)\sum_{e\in\cE_h^{ID}}\re\pdaj{\frac{\pa u_h^q}{\pa n_e}}{u_h^q}\nn\\
&\quad +2\sum_{e\in\cE_h^{ID}}\re\left(-\pd{\al\cdot n_e\av{\na u_h^q},\jm{\na u_h^q}}_e
+\pdaj{\frac{\pa u_h^q}{\pa n_e}}{v_h}\right) \nn\\
&\quad+2\sum_{e\in\cE_h^{ID}}\re\pdja{u_h^q}{\frac{\pa v_h}{\pa n_e}}
+2\im\Bigl(L_1(u_h^q,v_h)+\sum_{j=0}^q J_j(u_h^q,v_h)\Bigr).\nn
\end{align}

\medskip
{\em Step 2: Derivation of a reverse inequality.} Our task now is to estimate
each term on the right-hand side of \eqref{e3.8}. Since the terms on the
first four lines can be bounded in the exactly same way as done in
\cite{fw08a}, we omit their derivations and only give the final
results here for the reader's convenience.
\begin{align} \label{el1}
&2\re \bigl( (f,v_h)  +\langle g, v_h\rangle_{\Gamma_R} \bigr) \\
&\hskip 1in
\le C M(f,g)^2+\frac18\abs{u_h^q}_{1,h}^2
+\frac{c_{\Om_1}}{4}\sum_{e\in\cE_h^R}\norml{\na u_h^q}{e}^2.\nn \db\\
&2k^2\sum_{e\in\cE_h^I}\re\pd{\al\cdot n_e\av{u_h^q},\jm{u_h^q}}_e  \label{el2}\\
&\hskip 1in
\le \frac{k^2}{3}\norml{u_h^q}{\Om}^2
+C\sum_{e\in\cE_h^I}\frac{p\,k^2}{\ga_{0,e}}
\frac{\ga_{0,e}p}{h_e}\norml{\jm{u_h^q}}{e}^2. \nn \\
&k^2\pd{\al\cdot n_\Ome, |u_h^q|^2}_{\pa\Om}
\le C k^2 \norml{u_h^q}{\Ga_R}^2
+\sum_{e\in\cE_h^{D}}k^2\pd{\al\cdot n_e, |u_h^q|^2}_e. \label{el3}\db \\
&2k\im\pd{u_h^q,v_h}_{\Gamma_R}
-\sum_{e\in\cE_h^R}\pd{\al\cdot n_e, |\na u_h^q|^2}_e \label{el4} \\
&\hskip 1in
\le Ck^2\norml{u_h^q}{\Ga_R}^2
-\frac{c_{\Ome_1}}2 \sum_{e\in\cE_h^{R}}\norml{\na u_h^q}{e}^2. \nn\db\\
&2(d-1)\sum_{e\in\cE_h^{ID}}\re\pdaj{\frac{\pa u_h^q}{\pa n_e}}{u_h^q} \label{el5}\\
&\hskip 1in
\le \frac18\abs{u_h^q}_{1,h}^2+C\sum_{e\in\cE_h^{ID}}\frac{p}{\ga_{0,e}}\frac{\ga_{0,e}\,p}{h_e}\norml{\jm{u_h^q}}{e}^2.\nn \db
\end{align}

The extra work here is to estimate the  terms on the last two lines
in \eqref{e3.8}. For an edge/face $e\in\cE_h^{ID}$, let $\Om_e$  denote the set of element(s) in $\cT_h$ containing $e$ as one edge/face. By \eqref{eid3} we obtain
\begin{align}\label{e3.12}
2\sum_{e\in\cE_h^{ID}}&\re\left(-\pd{\al\cdot n_e\av{\na u_h^q},\jm{\na u_h^q}}_e
+\pdaj{\frac{\pa u_h^q}{\pa n_e}}{v_h}\right) \\
=&2\sum_{e\in\cE_h^{ID}}\sum_{\ell=1}^{d-1} \re
\int_e\left(\al\cdot\tau_e^\ell\av{\frac{\pa u_h^q}{\pa n_e}}
-\al\cdot n_e\av{\frac{\pa u_h^q}{\pa \tau_e^\ell}}\right)
\Jump{\frac{\pa\overline{u_h^q}}{\pa\tau_e^\ell}} \nn\\
\lesssim &\sum_{e\in\cE_h^{ID}}p^3h_e^{-\frac32}\sum_{K\in\Om_e}
\norml{\na u_h^q}{K}\norml{\jm{u_h^q}}{e}\nn\\
\le& \frac18\abs{u_h^q}_{1,h}^2+C\sum_{e\in\cE_h^{ID}}\frac{p^5}{\ga_{0,e}h_e^2}\frac{\ga_{0,e}\,p}{h_e}\norml{\jm{u_h^q}}{e}^2\nn
\end{align}
The first term on the line six of \eqref{e3.8} is bounded
as follows (compare with (4.20) of \cite{fw08a}):
\begin{align}\label{ID6}
2\sum_{e\in\cE_h^{ID}}\re\pdja{u_h^q}{\frac{\pa v_h}{\pa n_e}}
\lesssim& \sum_{e\in\cE_h^{ID}}p\,h_e^{-\frac12}\norml{\jm{u_h^q}}{e}\sum_{K\in\Om_e}\norml{\na v_h}{K} \\
\lesssim& \sum_{e\in\cE_h^{ID}}p^3h_e^{-\frac32}\norml{\jm{u_h^q}}{e}\sum_{K\in\Om_e}\norml{\na u_h^q}{K}\nn\\
\le &\frac18\abs{u_h^q}_{1,h}^2
+C\sum_{e\in\cE_h^{ID}}\frac{p^5}{\ga_{0,e}h_e^2}\frac{\ga_{0,e}\,p}{h_e}\norml{\jm{u_h^q}}{e}^2. \nn
\end{align}

The penalty term $L_1(\cdot,\cdot)$ is estimated as follows. Recall that $v_h|_K=\al\cdot\na u_h^q|_K$ with $\al = x-x_{\Om_1}$ for
each $K\in \cT_h$. Noting that
\begin{align}\label{ID1}
\frac{\pa v_h}{\pa \tau_e^\ell} & = \frac{\pa u_h^q}{\pa \tau_e^\ell}
+\al\cdot \nab \Bigl(\frac{\pa u_h^q}{\pa \tau_e^\ell}\Bigr)\\
&= \frac{\pa u_h^q}{\pa \tau_e^\ell} + \al\cdot n_e \frac{\p }{\p \tau_e^\ell}
\Bigl( \frac{\p u_h^q}{\p n_e} \Bigr)
+ \sum_{m=1}^{d-1} \al\cdot \tau_e^m \frac{\p }{\p \tau_e^m}
\Bigl( \frac{\p u_h^q}{\p \tau_e^\ell} \Bigr),  \quad 1\leq \ell \leq d-1,\nn\db
\end{align}
by the definition of $L_1(\cdot,\cdot)$ and Lemma~\ref{lem-inv} we get
\begin{align}\label{ID4}
&2\im L_1(u_h^q,v_h) = 2\im \sum_{e\in \cE_h^{ID}} \sum_{\ell=1}^{d-1}
\frac{\beta_{1,e} p}{h_e}\left(\pdjj{\frac{\pa u_h^q}{\pa \tau_e^\ell}}
{\al\cdot n_e \frac{\p}{\p \tau_e^\ell}\Bigl(\frac{\p u_h^q}{\p n_e}\Bigr)}\right. \\
&\hskip 1in \left.
+ \sum_{m=1}^{d-1} \pdjj{ \frac{\pa u_h^q}{\pa \tau_e^\ell}}
{\al\cdot\tau_e^m \frac{\p}{\p \tau_e^m}\left( \frac{\p u_h^q}{\p \tau_e^\ell} \right)}\right) \nn \db\\
&\leq 2\im \sum_{e\in \cE_h^D} \sum_{\ell=1}^{d-1}
\frac{\beta_{1,e} p}{h_e} \pd{\al\cdot n_e\,\frac{\pa u_h^q}{\pa \tau_e^\ell},
\frac{\p}{\p \tau_e^\ell} \Bigl( \frac{\p u_h^q}{\p n_e} \Bigr)}_e\nn\\
&\quad+C\sum_{e\in\cE_h^{I}} \sum_{\ell=1}^{d-1}
\frac{\beta_{1,e} p^3}{h_e^2} \norml{\jm{\frac{\pa u_h^q}{\pa \tau_e^\ell}}}{e}
\norml{\jm{\frac{\p u_h^q}{\p n_e}}}{e}\nn\\
&\quad+C\sum_{e\in\cE_h^{ID}} \sum_{\ell=1}^{d-1}
\frac{\beta_{1,e} p^3}{h_e^2} \norml{\jm{\frac{\pa u_h^q}{\pa \tau_e^\ell}}}{e}
\norml{\jm{\frac{\p u_h^q}{\p \tau_e^\ell}}}{e} \nn \db\\
&\leq 2\im \sum_{e\in \cE_h^D} \sum_{\ell=1}^{d-1}
\frac{\beta_{1,e} p}{h_e} \pd{\al\cdot n_e\,\frac{\pa u_h^q}{\pa \tau_e^\ell},
\frac{\p}{\p \tau_e^\ell} \Bigl( \frac{\p u_h^q}{\p n_e} \Bigr)}_e\nn\\
&\quad+C\sum_{e\in\cE_h^{I}} \sum_{\ell=1}^{d-1}
\frac{p^3}{h_e^2}\sqrt{\frac{\beta_{1,e}}{\ga_{1,e}}}
\left( \frac{\ga_{1,e} h_e}{p} \norml{\jm{\frac{\p u_h^q}{\p n_e}}}{e}^2
+ \frac{\beta_{1,e} p}{h_e}\norml{\jm{\frac{\p u_h^q}{\p \tau_e^\ell}}}{e}^2\right)\nn\\
&\quad+C\sum_{e\in\cE_h^{ID}} \sum_{\ell=1}^{d-1}
\frac{p^2}{h_e} \frac{\beta_{1,e} p}{h_e} \norml{\jm{\frac{\p u_h^q}{\p \tau_e^\ell}}}{e}^2
.\nn
\end{align}
We remark that $\im L_1(u_h^q,v_h)=0$ when $p=q=1$.

Next we estimate the penalty terms $ J_j(u_h^q,v_h)$. Since $u_h^q$ and $v_h$ are piecewise polynomials of
degree $p$ in general and those terms contain jumps of high order
normal derivatives, it is quite delicate to control those terms
as shown below.

By direct calculations we get that on each
edge/face $e$ of $K\in \cT_h$
\begin{align}
\frac{\pa^j v_h}{\pa n_e^j} & = j\frac{\pa^j u_h^q}{\pa n_e^j}
+ \al \cdot \nab \Bigl(\frac{\pa^j u_h^q}{\pa n_e^j}\Bigr) \label{ID2} \\
& =j \frac{\pa^j u_h^q}{\pa n_e^j}+\al\cdot n_e \frac{\pa^{j+1} u_h^q}{\pa n_e^{j+1}}
+ \sum_{m=1}^{d-1}   \al\cdot \tau_e^m \frac{\p}{\p \tau_e^m}
\Bigl(\frac{\p^j u_h^q}{\p n_e^j} \Bigr), \quad1\leq j\leq p-1,\nn \db\\
\frac{\pa^p v_h}{\pa n_e^p} & = p\frac{\pa^p u_h^q}{\pa n_e^p}
+ \al\cdot \nab \Bigl(\frac{\pa^p u_h^q}{\pa n_e^p}\Bigr)
= p\frac{\pa^p u_h^q}{\pa n_e^p}. \label{ID3}
\end{align}
Here we have used the fact that $(p+1)$th order derivatives of $u_h^q$ is
zero because $u_h^q$ is a polynomial of degree at most $p$.

For $j=1,2,\cdots q-1$, by \eqref{ID2} we have
\begin{align} \label{ID5}
2\im J_j(u_h^q,v_h)&=2\im\sum_{e\in\cE_h^I}\ga_{j,e}
\left(\frac{h_e}{p}\right)^{2j-1} \left(
\al\cdot n_e \pdjj{\frac{\p^j u_h^q}{\p n_e^j}}{\frac{\p^{j+1} u_h^q}{\p n_e^{j+1}}}
\right. \\
&\hskip 0.7in \left.
+\sum_{\ell=1}^{d-1} \al\cdot \tau_e^\ell \pdjj{\frac{\p^j u_h^q}{\p n_e^j}}
{\frac{\p }{\p \tau_e^\ell} \Bigl( \frac{\p^j u_h^q}{\p n_e^j} \Bigr)} \right) \nn\db \\
&\lesssim \sum_{e\in\cE_h^I} \gamma_{j,e} \left(\frac{h_e}{p}\right)^{2j-1}
\norml{\jm{\frac{\p^j u_h^q}{\p n_e^j}}}{e}\nn\\
&\hskip 0.7in  \times\left(\norml{\jm{\frac{\p^{j+1} u_h^q}{\p n_e^{j+1}}}}{e}
+\frac{p^2}{h_e}\norml{\jm{\frac{\p^j u_h^q}{\p n_e^j}}}{e} \right) \nn\db\\
&\lesssim \sum_{e\in\cE_h^I} \frac{p}{h_e}\sqrt{\frac{\ga_{j,e}}{\ga_{j+1,e}}}
\left(\gamma_{j,e} \left(\frac{h_e}{p}\right)^{2j-1}\norml{\jm{\frac{\p^j u_h^q}{\p n_e^j}}}{e}^2\right.  \nn  \\
&\hskip 0.7in
\left.+\gamma_{j+1,e} \left(\frac{h_e}{p}\right)^{2j+1}\norml{\jm{\frac{\p^{j+1} u_h^q}{\p n_e^{j+1}}}}{e}^2\right)\nn\\
&\hskip 0.7in+\sum_{e\in\cE_h^I} \frac{p^2}{h_e}\gamma_{j,e} \left(\frac{h_e}{p}\right)^{2j-1}\norml{\jm{\frac{\p^j u_h^q}{\p n_e^j}}}{e}^2. \nn
\end{align}
If $q<p$, then, from Lemma~\ref{lem-inv} and the inequality 
$\norml{\dfrac{\pa^q \vp}{\pa n_e^q}}{\pa K}
\lesssim ph_e^{-\frac12}\abs{\vp}_{H^q(K)}$, we have
\begin{align} \label{ID5q}
\im 2J_q(u_h^q,v_h) &\lesssim\sum_{e\in\cE_h^I}\ga_{q,e} \left(\frac{h_e}{p}\right)^{2q-1}
\norml{\frac{\p^q u_h^q}{\p n_e^q}}{e}p\,h_e^{-\frac12}\sum_{K\in\Om_e}\abs{v_h}_{H^q(K)}  \\
&\lesssim\sum_{e\in\cE_h^I}\ga_{q,e}\,\left(\frac{h_e}{p}\right)^{2q-1}\,\frac{p^{2q+1}}{h_e^{q+\frac12}} \norml{\frac{\p^q u_h^q}{\p n_e^q}}{e}\sum_{K\in\Om_e}\norml{\na u_h^q}{K}\nn  \\
&\le\frac18\abs{u_h^q}_{1,h}^2+C\sum_{e\in\cE_h^I}\frac{\ga_{q,e}\,p^{2q+3}}{h_e^2}\ga_{q,e} \left(\frac{h_e}{p}\right)^{2q-1} \norml{\frac{\p^q u_h^q}{\p n_e^q}}{e}^2. \nn
\end{align}
If $q=p$, \eqref{ID3} and the definition of $J_p(\cdot,\cdot)$ immediately
imply that (compare with (4.14) of \cite{fw08a})
\begin{equation}\label{e3.7a}
2\im J_p(u_h^q,v_h)=2\im J_p(u_h^q,u_h^q) =0.
\end{equation}
The estimate for $\im J_0(u_h^q,v_h)$ is similar to \eqref{ID5}, so we get
\begin{align} \label{ID5a}
&2\im J_0(u_h^q,v_h)
=2\im\sum_{e\in\cE_h^{ID}}\frac{\ga_{0,e}\,p}{h_e}\pdjj{u_h^q}{\al\cdot n_e\frac{\pa u_h^q}{\pa n_e}
+\sum_{j=1}^{d-1}\al\cdot\tau_e^j\frac{\pa u_h^q}{\pa \tau_e^j}} \db \\
&\le 2\im\sum_{e\in\cE_h^D}\frac{\ga_{0,e}\,p}{h_e}\pd{\al\cdot n_e u_h^q,\frac{\pa u_h^q}{\pa n_e}}_e\nn\\
&\qquad+C\sum_{e\in\cE_h^I}\frac{\ga_{0,e}\,p}{h_e} \norml{\jm{u_h^q}}{e}
\norml{\jm{\frac{\pa u_h^q}{\pa n_e}}}{e} \nn\\
&\qquad+C\sum_{e\in\cE_h^{ID}}\frac{\ga_{0,e}\,p}{h_e}
\norml{\jm{u_h^q}}{e}\frac{p^2}{h_e}\norml{\jm{u_h^q}}{e}
\nn\db\\
&\le 2\im\sum_{e\in\cE_h^D}\frac{\ga_{0,e}\,p}{h_e}\pd{\al\cdot n_e u_h^q,\frac{\pa u_h^q}{\pa n_e}}_e
+C\sum_{e\in\cE_h^{ID}}\frac{p^2}{h_e}\frac{\ga_{0,e}\,p}{h_e}
\norml{\jm{u_h^q}}{e}^2 \nn\\
&\qquad+C\sum_{e\in\cE_h^I}\frac{p}{h_e}\sqrt{\frac{\ga_{0,e}}{\ga_{1,e}}}
\left(\frac{\ga_{0,e}\,p}{h_e}\norml{\jm{u_h^q}}{e}^2+
\frac{\ga_{1,e}h_e}{p}\norml{\Jump{\frac{\pa u_h^q}{\pa n_e}}}{e}^2\right).\nn
\end{align}

We also need the following estimate (compare with (4.22) of \cite{fw08a})
\begin{align}
&\sum_{e\in\cE_h^{D}} \bigg(k^2\pd{\al\cdot n_e, |u_h^q|^2}_e
+\pd{\al\cdot n_e, |\na u_h^q|^2}_e\label{el6}\\
& \qquad+ \sum_{\ell=1}^{d-1}
\frac{2\beta_{1,e} p}{h_e} \im \pd{\al\cdot n_e\,\frac{\pa u_h^q}{\pa \tau_e^\ell},
\frac{\p}{\p \tau_e^\ell} \Bigl( \frac{\p u_h^q}{\p n_e} \Bigr)}_e \nn\\
& \qquad+\frac{2\ga_{0,e}\,p}{h_e}\im\pd{\al\cdot n_e u_h^q,\frac{\pa u_h^q}{\pa n_e}}_e \bigg) \nn \db\\
&\le -\sum_{e\in\cE_h^{D}}\bigg\langle\al\cdot n_D, k^2|u_h^q|^2+|\na u_h^q|^2-2\sum_{\ell=1}^{d-1}
\frac{\beta_{1,e} p}{h_e} \Big|\frac{\pa u_h^q}{\pa \tau_e^\ell}\Big|\Big|\frac{\p}{\p \tau_e^\ell} \Bigl( \frac{\p u_h^q}{\p n_e} \Bigr)\Big|\nn\\
&\qquad-2\frac{\ga_{0,e}\,p}{h_e}|u_h^q||\na u_h^q|\bigg\rangle_e \nn \db\\
&\leq  -c_D\sum_{e\in\cE_h^{D}} \Bigl( k^2\|u_h^q\|_{L^2(e)}^2
+ \frac12 \|\na u_h^q\|_{L^2(e)}^2\Bigr) \nn \\
&\quad
+C\sum_{e\in\cE_h^{D}}\frac{\beta_{1,e}\,p^5}{h_e^3}\sum_{\ell=1}^{d-1} \frac{\beta_{1,e}\,p}{h_e}\norml{\frac{\pa u_h^q}{\pa \tau_e^\ell}}{e}^2+C\sum_{e\in\cE_h^{D}}\frac{\ga_{0,e}\,p}{h_e} \frac{\ga_{0,e}\,p}{h_e}\norml{u_h^q}{e}^2,\nn
\end{align}
where we have used the inverse inequality and the assumption that $D$ is star-shaped to derive the last inequality.

\medskip
{\em Step 3: Finishing up.} We only prove the case of $q=p$ since the
proof for $q<p$ is the same except using \eqref{ID5q} instead of \eqref{e3.7a}.
Substituting \eqref{ID5}, \eqref{e3.7a}, \eqref{ID5a} (with $q=p$) and
 \eqref{el1}--\eqref{el5} into \eqref{e3.8}, and using \eqref{el6} we obtain
\begin{align*}
2k^2&\norml{u_h^p}{\Om}^2\\
\le& (d-2)\re\bigl( (f,u_h^p) +\pd{g, u_h^p}_{\Ga_R} \bigr)
+C M(f,g)^2+ \frac{k^2}{3}\norml{u_h^p}{\Om}^2\db\\
&-\frac38\abs{u_h^p}_{1,h}^2+\abs{u_h^p}_{1,h}^2-2\sum_{e\in\cE_h^{ID}}\re\pdaj{\frac{\pa u_h^p}{\pa n_e}}{u_h^p}\db\\
&-\frac{c_{\Ome_1}}4 \sum_{e\in\cE_h^R}\norml{\na u_h^p}{e}^2
+ C k^2 \norml{u_h^p}{\Ga_R}^2 \db\\
&-c_D\sum_{e\in\cE_h^{D}} \Bigl( k^2\|u_h^p\|_{L^2(e)}^2
+ \frac12 \|\na u_h^p\|_{L^2(e)}^2\Bigr)\db\\
&+C\sum_{e\in\cE_h^I}\left(\frac{p\,(k^2+1)}{\ga_{0,e}} +\frac{p^5}{\ga_{0,e}h_e^2}
+\frac{p^2}{h_e}\right)\frac{\ga_{0,e}\,p}{h_e}\norml{\jm{u_h^p}}{e}^2 \db\\
&+C\sum_{e\in\cE_h^D}\left(
\frac{\ga_{0,e}\,p}{h_e}
+\frac{p^5}{\ga_{0,e}h_e^2}
+\frac{p^2}{h_e}\right)\frac{\ga_{0,e}\,p}{h_e}\norml{\jm{u_h^p}}{e}^2\db\\
&+C\sum_{e\in\cE_h^{D}}
\frac{\beta_{1,e}\,p^5}{h_e^3}\sum_{\ell=1}^{d-1} \frac{\beta_{1,e}\,p}{h_e}\norml{\frac{\pa u_h^p}{\pa \tau_e^\ell}}{e}^2\db\\
&+C\sum_{e\in\cE_h^{I}} \sum_{\ell=1}^{d-1}
\frac{p^3}{h_e^2}\sqrt{\frac{\beta_{1,e}}{\ga_{1,e}}}
\left( \frac{\ga_{1,e} h_e}{p} \norml{\jm{\frac{\p u_h^p}{\p n_e}}}{e}^2
+ \frac{\beta_{1,e} p}{h_e}\norml{\jm{\frac{\p u_h^p}{\p \tau_e^\ell}}}{e}^2\right)\db\\
&+C\sum_{e\in\cE_h^{ID}} \sum_{\ell=1}^{d-1}
\frac{p^2}{h_e} \frac{\beta_{1,e} p}{h_e} \norml{\jm{\frac{\p u_h^p}{\p \tau_e^\ell}}}{e}^2\db\\
&+C\sum_{j=1}^{p-1}\sum_{e\in\cE_h^I} \frac{p^2}{h_e}\gamma_{j,e} \left(\frac{h_e}{p}\right)^{2j-1}\norml{\jm{\frac{\p^j u_h^p}{\p n_e^j}}}{e}^2 \db\\
&+C\sum_{j=0}^{p-1}\sum_{e\in\cE_h^{I}} \frac{p}{h_e}\sqrt{\frac{\ga_{j,e}}{\ga_{j+1,e}}}\left(
\gamma_{j,e} \left(\frac{h_e}{p}\right)^{2j-1}\norml{\jm{\frac{\p^j u_h^p}{\p n_e^j}}}{e}^2\right.\\
&\left.\hskip 1.5in+
\gamma_{j+1,e} \left(\frac{h_e}{p}\right)^{2j+1}\norml{\jm{\frac{\p^{j+1} u_h^p}{\p n_e^{j+1}}}}{e}^2\right). \nn
\end{align*}
Therefore, it follows from Lemma~\ref{lem3.1} and \eqref{csta} that
\begin{align*}
&2k^2\norml{u_h^p}{\Om}^2+\frac38\abs{u_h^p}_{1,h}^2
+\frac{c_{\Ome_1}}4 \sum_{e\in\cE_h^R}\norml{\na u_h^p}{e}^2 \\
&\qquad + c_D\sum_{e\in\cE_h^{D}} \Bigl( k^2\|u_h^p\|_{L^2(e)}^2
+ \frac12 \|\na u_h^p\|_{L^2(e)}^2\Bigr) \db\\
&\le C M(f,g)^2+\frac{4k^2}{3}\norml{u_h^p}{\Om}^2
+Ck^2\csta{p} \bigl|(f,u_h^p)+\langle g, u_h^p\rangle_{\Ga_R}\bigr|\\
&\le C k^2\csta{p}^2 M(f,g)^2+\frac{5k^2}{3}\norml{u_h^p}{\Om}^2,
\end{align*}
where we have used the following inequality, which is a consequence
of \eqref{e3.2},
\begin{equation*}\label{3.6t}
k^2\norml{u_h^p}{\Ga_R}^2\le k^2\norml{u_h^p}{\Om}^2+M(f,g)^2
\end{equation*}
to derive the last inequality. $M(f,g)$ and $\csta{p}$ are defined
by \eqref{Mfg} and \eqref{csta}, respectively. Hence,
\begin{align*}
\norml{u_h^p}{\Om} &+ \frac{1}{k} \abs{u_h^p}_{1,h}
+\frac{1}{k} \biggl( c_{\Om_1}\sum_{e\in\cE_h^R}\norml{\na u_h^p}{e}^2 \biggr)^{\frac12} \\
&+ \frac{1}{k} \biggl( \sum_{e\in\cE_h^{D}} c_D\Bigl( k^2\|u_h^p\|_{L^2(e)}^2
+  \|\na u_h^p\|_{L^2(e)}^2\Bigr)  \biggr)^{\frac12}
\lesssim \csta{p} M(f,g),
\end{align*}
which together with \eqref{e3.2} gives \eqref{e3.3}.  The proof is completed.
\end{proof}

As \eqref{edg} can be written as a linear system,
an immediate consequence of the above stability estimate
is the following well-posedness theorem for \eqref{edg}.

\begin{theorem}\label{existence}
For $k>0$ and $h>0$, the $hp$-IPDG method \eqref{edg}  has a
unique solution $u_h^q$ provided that $\ga_{0,e}, \ga_{1,e},\cdots \ga_{q,e}> 0$ and $\beta_{1,e}\ge 0$.
\end{theorem}

Next we consider the case of quasi-uniform meshes.
Note that large  penalty parameters ($\ga_j, j\ge 1$)
for jumps of normal derivatives may cause a large interpolation
error in the norm $\norm{\cdot}_{1,h,q}$, and hence, may pollute the
error estimates of the IPDG solution (see Section~\ref{sec-err}).
It is interested to minimize the stability constant $\csta{q}$
under the constraints of $\beta_{1,e}\ge 0$ and
$\frac{p}{\ga_{0,e}}+\sum_{j=1}^qp^{2j-1}\ga_{j,e}\lesssim 1.$ We have the following
consequence of Theorem~\ref{thm_sta}.   The proof  is
straightforward and is omitted.

\begin{theorem}\label{thm_sta_uniform}
Let $h=\max h_e.$ Suppose the mesh $\cT_h$
is quasi-uniform, that is $h_e\simeq h.$ Suppose $k\gtrsim 1$
and $k\,h\lesssim 1$.  Assume that $\ga_{j,e}\simeq\ga_j,  j=0,1,\cdots,q$,
$\sum_{j=1}^qp^{2j-1}\ga_j\lesssim 1$, $0\le\beta_{1,e}\lesssim\frac{h^2}{p^4}\ga_0$, and that
\[
\left\{\begin{array}{ll}
\ga_0\simeq  p^{\frac73}h^{-\frac{2}{3}},\quad \ga_j\gtrsim p^{-\frac{10}{3}}h^{\frac{2}{3}}\,\ga_{j-1} & \text{if } q=p,\\
\ga_0\simeq \min\set{p^{\frac{3q+1}{q+1}}h^{-\frac{q}{q+1}},p^{\frac73}h^{-\frac{2}{3}}},
\quad \ga_j\gtrsim \left(\frac{\ga_0\,h}{p^4}\right)^2\ga_{j-1}, \quad \ga_q\lesssim \frac{1}{\ga_0p^{2q-2}} & \text{if } q< p.
\end{array}\right.
\]
Suppose $D=\emptyset$ if $q\ge 2$ or $q=p=1$. Then
\begin{equation*}
\norml{u_h}{\Om}+\frac1k\norm{u_h}_{1,h,q}\lesssim\csta{q} M(f,g),
\end{equation*}
where
\[
\csta{q}\lesssim\left\{\begin{aligned}
&\frac{p^{\frac{8}{3}}}{k^2 h^{\frac43}} &\text{ if } q= p,\\
&\max\set{\frac{p^{\frac{8}{3}}}{k^2 h^{\frac43}},\frac{p^{\frac{2q+4}{q+1}}}{k^2 h^{\frac{q+2}{q+1}}}} &\text{ if } q<p.
\end{aligned}
\right.
\]
\end{theorem}

It is clear that, in the above theorem,
$\ga_0\simeq p^{\frac73}h^{-\frac{2}{3}}$ if $q=p$ or $q\ge 2$,
$\ga_0\simeq p^{\frac{3q+1}{q+1}}h^{-\frac{q}{q+1}}$ otherwise.
We conclude this section by several remarks.

\begin{remark}
(a) The $hp$-IPDG method \eqref{edg} is well-posed for all $h,k>0$
provided that all penalty parameters are positive. As a comparison, 
we recall that the standard finite element method is well-posed only 
if mesh size $h$ satisfies a constraint $h=O(k^{-\rho})$ for 
some $\rho\geq 1$, hence, the existence is only guaranteed for very 
small mesh size $h$ when wave number $k$ is large.

(b). It is well known that \cite{arnold82,abcm01,baker77,fk07,rwg99,w78}
symmetric IPDG methods for coercive elliptic and parabolic PDEs
often require the penalty parameter $\gamma_{0,e}$ is sufficiently
large to guarantee the well-posedness of numerical solutions, and the
low bound for $\gamma_{0,e}$ is theoretically hard to determine and
is also problem-dependent. However, this is no issue for
scheme \eqref{edg}, which solves the (indefinite) Helmholtz
equation, because they are well-posed for all $\gamma_{0,e}>0$.

(c) In the linear element case of $q=p=1$, a better estimate is obtained in \cite{fw08a} with the help of the penalty term $L_1$. Unfortunately the penalty term $L_1$ does not help too much for higher order elements.

(d). The stability estimates will be improved greatly when  $k^3h^2p^{-1}\le C_0$, where $C_0$ is some constant independent of $k$, $h$, $p$, and the penalty parameters (see Theorem~\ref{thm_main2dg0}).  The above estimates
are only interested  when $k\,h\lesssim 1$ and $k^3 h^2p^{-1}\gtrsim 1$.
\end{remark}

%%%%%%%%%%%%%%%%%%%%%%%%
\section{Error estimates}\label{sec-err}
In this section, we derive the error estimates of the solutions
of scheme \eqref{edg}. This will be done in two steps.
First, we introduce  elliptic projections of the PDE solution
$u$ and derive error estimates for the projections.
We note that such a result also has an independent interest.
Second, we bound the error between the projections and the IPDG
solutions by making use of the stability results obtained in
Section~\ref{sec-sta}. In this and the next section, we assume that the
mesh $\cT_h$ is quasi-uniform, that $\ga_{j,e}\simeq\ga_j> 0$ for $j=0,1,\cdots,q$,
and that $\beta_{1,e}\simeq \beta_1\ge 0$. Let $h=\max h_e$. For simplicity, we also assume that the
mesh $\cT_h$ is  conforming, that is, $\cT_h$ contains no hanging nodes, since the parallel results for non-conforming meshes can be derived in a similar way.

In this and the next section  the following assumption (cf. Theorem~\ref{stability}) is required for some results.
Suppose problem \eqref{e1.1}--\eqref{e1.3} is $H^{s}$-regular and
\begin{equation}\label{Areg}
\norm{u}_{H^{s}(\Om)}\lesssim k^{s-1} M_s(f,g),
\end{equation}
where $s\ge 2$ is an integer and
\[
M_s(f,g)= \norm{f}_{H^{\max\{0,s-2\}}(\Ome)} + \|g\|_{L^2(\Ga_R)}.
\]
Note that the assumption is proved to hold for $s=2$ and $M_2(f,g)=M(f,g)$.

%%%%%%%%%%%%%
\subsection{Elliptic projection and its error estimates}
For any $w\in E^q\cap H^1_{\Ga_D}(\Ome)\cap H^2_{\mbox{\tiny loc}} (\Ome)$,
we define its elliptic projection $\tilde{w}_h^q\in V_h^p$ by
\begin{equation}\label{e4.3}
a_h^q(\tilde{w}_h^q,v_h)+\i k\pd{\tilde{w}_h^q,v_h}_{\Gamma_R}
=a_h^q(w,v_h)+\i k\pd{w,v_h}_{\Gamma_R}
\qquad\forall v_h\in V_h^p.
\end{equation}
In other words, $\tilde{w}_h$ is an IPDG approximation to
the solution $w$ of the following (complex-valued) Poisson problem:
\begin{alignat*}{2}
-\Del w &= F &&\qquad \mbox{in } \Ome,\\
\frac{\p w}{\p n_{\Gamma_R}} +\i k w &=\psi &&\qquad \mbox{on }\Gamma_R,\\
w &= 0&&\qquad \mbox{on }\Gamma_D
\end{alignat*}
for some given functions $F$ and $\psi$ which are determined by $w$.

Before estimating the projection error, we state the following continuity and
coercivity properties for the sesquilinear form $a_h^q(\cdot,\cdot)$.
Since they follow easily from \eqref{eah}--\eqref{ea2}, so
we omit their proofs to save space.
\begin{lemma}\label{lem4.1}
For any $v\in E^q$ and $w\in E^q\cap H^1_{\Ga_D}$, the mesh-dependent
sesquilinear form $a_h^q(\cdot,\cdot)$ satisfies
\begin{equation}
|a_h^q(v,w)|, |a_h^q(w,v)|\lesssim \norm{v}_{1,h,q} \normeq{w}. \label{e4.1}
\end{equation}
In addition, for any $0<\ep<1$, there exists a positive constant
$c_\ep$ independent of $k$, $h$, $p$, and the penalty parameters such that
\begin{equation}
\re a_h^q(v_h,v_h)+\Bigl(1-\ep+\frac{c_\ep\,p}{\ga_0}\Bigr)
\im a_h^q(v_h,v_h)\ge (1-\ep)\norm{v_h}_{1,h,q}^2\qquad
\forall v_h\in V_h^p. \label{e4.2}
\end{equation}
\end{lemma}

Let $u$ be the solution of problem \eqref{e1.1}--\eqref{e1.3}
and $\tuh^q$ be its elliptic projection defined as above.
Define $\eta:=u-\tuh^q$. Then \eqref{e4.3} immediately implies the following
Galerkin orthogonality:
\begin{equation}\label{e4.4}
a_h^q(\eta,v_h)+\i k\pd{\eta,v_h}_{\Gamma_R} =0 \qquad\forall v_h\in V_h^p.
\end{equation}

\begin{lemma}\label{lem-4.2a}
Suppose problem \eqref{e1.1}--\eqref{e1.3} is $H^{\max\set{q+1,2}}$-regular.
Then there hold the following estimates:
\begin{align}
\norm{u-\tuh^q}_{1,h,q}& + \sqrt{\lambda k} \norml{u-\tuh^q}{\Gamma_R}\label{e4.2aa}\\
\lesssim&\inf_{z_h\in V_h^p\cap H_{\Ga_D}^1(\Om)} \left(\lambda \normeq{u-z_h}
+\sqrt{\lambda k}\norml{u-z_h}{\Gamma_R}\right),\nn\\
\norml{u-\tuh^q}{\Om}
\lesssim & h\,\Big(1+\frac{1}{\ga_0\,p}+\frac{\gamma_1}{p}+kh\Big)^{\frac12}\label{e4.2ab}\\
&\times\inf_{z_h\in V_h^p\cap H_{\Ga_D}^1(\Om)} \left(\lambda \normeq{u-z_h}
+\sqrt{\lambda k}\norml{u-z_h}{\Gamma_R}\right),\nn
\end{align}
where $\lambda:=1+\frac{p}{\gamma_0}$ and $\ga_1=0$ if $q=0$.
\end{lemma}
\begin{proof} For any $z_h\in V_h^p\cap H_{\Ga_D}^1(\Om)$, let $\eta_h=\tuh^q-z_h$. From $\eta_h+\eta=u-z_h$ and \eqref{e4.4}, we have
\begin{equation}\label{e4.5}
a_h^q(\eta_h,\eta_h)+\i k\pd{\eta_h,\eta_h}_{\Gamma_R}
=a_h^q(u-z_h,\eta_h)+\i k\pd{u-z_h,\eta_h}_{\Gamma_R}.
\end{equation}
Take  $\ep=\frac12$ in \eqref{e4.2} and assume without loss of generality
that $c_{\frac12}>\frac12$. It follows from \eqref{e4.2} and \eqref{e4.5} that
\begin{align*}
\frac12\norm{\eta_h}_{1,h,q}^2\le& \re a_h^q(\eta_h,\eta_h)
+\Bigl(\frac12+\frac{c_{\frac12}\,p}{\ga_0}\Bigr)\im a_h^q(\eta_h,\eta_h)\\
=&\re\bigl(a_h^q(u-z_h,\eta_h)+\i k\pd{u-z_h,\eta_h}_{\Gamma_R}\bigr)
-\Bigl(\frac12+\frac{c_{\frac12}\,p}{\ga_0}\Bigr)k\pd{\eta_h,\eta_h}_{\Gamma_R} \\
&\qquad +\Bigl(\frac12+\frac{c_{\frac12}\,p}{\ga_0}\Bigr)
\im\left(a_h^q(u-z_h,\eta_h)+\i k\pd{u-z_h,\eta_h}_{\Gamma_R}\right)\\
\le &C\lambda \Bigl(\norm{\eta_h}_{1,h,q}\normeq{u-z_h}
+k\norml{u-z_h}{\Gamma_R}^2\Bigr)
-\frac{\lambda k}{4}\norml{\eta_h}{\Gamma_R}^2.
\end{align*}
Therefore,
\begin{align}\label{e4.6a}
\norm{\eta_h}_{1,h,q}^2 +\lambda k\norml{\eta_h}{\Gamma_R}^2
\lesssim &\lambda^2 \normeq{u-z_h}^2
+\lambda k\norml{u-z_h}{\Gamma_R}^2
\end{align}
which together with $\eta=u-z_h-\eta_h$ yields \eqref{e4.2aa}.

To show \eqref{e4.2ab}, we use the Nitsche's duality
argument (cf. \cite{bs94,ciarlet78}). Consider the following auxiliary problem:
\begin{alignat}{2}\label{e4.7}
-\De w &=\eta &&\qquad\text{in }\Om,\\
\frac{\pa w}{\pa n}-\i k w &=0 &&\qquad\text{on }\Gamma_R, \nn \\
w &=0 &&\qquad\text{on }\Gamma_D. \nn
\end{alignat}
It can be shown that $w$ satisfies
\begin{equation}\label{e4.8}
\abs{w}_{H^2(\Om)}\lesssim \norml{\eta}{\Om}.
\end{equation}
Let $\hat{w}_h^1$ be the continuous linear finite element interpolant
of $w$ on $\cT_h$. From \eqref{e4.4},
\begin{equation*}
   a_h^1(\eta,\hat{w}_h^1)+\i k\pd{\eta,\hat{w}_h^1}_{\Gamma_R}= a_h^q(\eta,\hat{w}_h^1)+\i k\pd{\eta,\hat{w}_h^1}_{\Gamma_R} =0.
\end{equation*}
Testing the conjugated
\eqref{e4.7} by $\eta$ and using the above orthogonality we get
\begin{align*}
\norml{\eta}{\Om}^2 &=a_h^1(\eta,w)+\i k\pd{\eta,w}_{\Gamma_R} \\
&=a_h^1(\eta,w-\hat{w}_h^1)+\i k\pd{\eta,w-\hat{w}_h^1}_{\Gamma_R} \\
&\lesssim \norm{\eta}_{1,h,1}\normeone{w-\hat{w}_h^1}
+k\norml{\eta}{\Gamma_R}\norml{w-\hat{w}_h^1}{\Gamma_R}. \\
&\lesssim \norm{\eta}_{1,h,q}
\Big(1+\frac{1}{\ga_0\,p}+\frac{\gamma_1}{p}\Big)^{\frac12}\,h\abs{w}_{H^2(\Om)}
+k\norml{\eta}{\Gamma_R}\,h^{\frac32}\abs{w}_{H^2(\Om)},
\end{align*}
which together with \eqref{e4.2aa} and \eqref{e4.8} gives \eqref{e4.2ab}.
The proof is completed.
\end{proof}

We have the following lemma that gives approximation properties of the space $V_h^p\cap H_{\Ga_D}^1(\Om)$.
\begin{lemma}\label{lem-app2}$\,$
\begin{enumerate}
  \item[(i)]Let $\mu= \min\set{p+1,s}$ and $q<\mu$. Suppose $u\in H^s(\Om)\cap H_{\Ga_D}^1(\Om)$. Then there exists $\hat u_h\in V_h^p\cap H_{\Ga_D}^1(\Om)$ such that
\begin{align}
&\norml{u-\hat{u}_h}{\Gamma_R}
\lesssim\frac{h^{\mu-\frac12}}{p^{s-\frac12}}\norm{u}_{H^{s}(\Om)},
\label{e4.6c}\\
&\normeq{u-\hat u_h}\lesssim \Big(1+\frac{p}{\ga_0}+\sum_{j=1}^qp^{2j-1}\ga_j\Big)^{\frac12}\;\frac{h^{\mu-1}}{p^{s-1}}\norm{u}_{H^{s}(\Om)}. \label{e4.6b}
\end{align}
  \item[(ii)] Suppose $u\in H^{\max\set{q+1,2}}(\Om)\cap H_{\Ga_D}^1(\Om)$. Then there exists $\hat u_h\in V_h^p\cap H_{\Ga_D}^1(\Om)$ such that \eqref{e4.6c} holds with $s=2$ and
\begin{align}
&\normeq{u-\hat u_h}\lesssim \Big(1+\frac{p}{\ga_0}+ p\,\ga_1+\sum_{j=2}^q p^{2j-2}\ga_j\Big)^{\frac12}\;\frac{h}{p}\norm{u}_{H^{2}(\Om)}, \label{e4.6d}
\end{align}
\end{enumerate}
where $\ga_1=0$ if $q=0.$
\end{lemma}
\begin{proof}
The following $hp$ approximation properties are well-known for the $hp$
finite element functions (cf. \cite{bs87,guo06,gs07}):
\begin{itemize}
\item There exists $\check u_h\in V_h^p$ such that, for $ j=0,1,\cdots,s$,
\begin{equation}\label{elem-app-a}
\norm{u-\check u_h}_{H^j(\cT_h)}
:=\Big(\sum_{K\in\cT_h}\norm{u-\check u_h}_{H^j(K)}^2\Big)^{\frac12}\lesssim\frac{h^{\mu-j}}{p^{s-j}}\norm{u}_{H^{s}(\Om)}.
 \end{equation}
\item There exists $\hat u_h\in V_h^p\cap H_{\Ga_D}^1(\Om)$ such that
 \begin{equation}\label{elem-app2-a}
   \norm{u-\hat u_h}_{H^j(\Om)}\lesssim \frac{h^{\mu-j}}{p^{s-j}}\norm{u}_{H^{s}(\Om)}, \quad j=0,1.
 \end{equation}
\end{itemize}
Here the invisible constants in the above two inequalities depending on $s$ but independent of $k$, $h$, $p$, and the penalty parameters. Then \eqref{e4.6c} follows from \eqref{elem-app2-a} and the trace inequality.

It follows from the inverse inequalities in Lemma~\ref{lem-inv} that, for $1\le j\le q+1$,
\begin{align}\label{elem-app2-c}
\norm{u-\hat u_h}_{H^j(\cT_h)}
&\le \norm{u-\check u_h}_{H^j(\cT_h)}+\norm{\check u_h-\hat u_h}_{H^j(\cT_h)}\\
&\lesssim \norm{u-\check u_h}_{H^j(\cT_h)}
+\frac{p^{2(j-1)}}{h^{j-1}}\norm{\check u_h-\hat u_h}_{H^1(\Om)} \nn \\
&\lesssim \frac{h^{\mu-j}}{p^{s-2j+1}}\norm{u}_{H^{s}(\Om)}.\nn
\end{align}
Therefore, by the following local trace inequality
\begin{equation*}\label{elti}
\norml{v}{\pa K}^2\lesssim h_K^{-1}\norml{v}{K}^2+\norml{v}{K}\norml{\na v}{K},
\end{equation*}
we have
\begin{align}\label{elem-app2-b}
\sum_{K\in\cT_h}&\sum_{e\subset\pa K}\norml{\frac{\pa^j (u-\hat u_h)}{\pa n_e^j}}{e}^2\\
\lesssim&h^{-1}\norm{u-\hat u_h}_{H^j(\cT_h)}^2
+\norm{u-\hat u_h}_{H^j(\cT_h)}\norm{u-\hat u_h}_{H^{j+1}(\cT_h)}\nn\\
\lesssim&\frac{h^{2\mu-2j-1}}{p^{2s-4j}}\norm{u}_{H^{s}(\Om)}^2.\nn
\end{align}
Noting that $\hat u_h$ is continuous, we have from
\eqref{elem-app2-a} and \eqref{elem-app2-b},
\begin{align}\label{elem-app2-h}
\normeq{u-\hat u_h}^2=&\abs{u-\hat u_h}_{1,h}^2
+ \sum_{j=1}^q\sum_{e\in\cE_h^I} \ga_{j,e} \left(\frac{h_e}{p}\right)^{2j-1}
\norml{\jm{\frac{\pa^j (u-\hat u_h)}{\pa n_e^j}}}{e}^2\\
&+\sum_{e\in\cE_h^{ID}} \frac{h_e}{\ga_{0,e}\,p}\norml{\av{\frac{\pa (u-\hat u_h)}{\pa n_e}}}{e}^2\nn \\
\lesssim &\Big(1+\frac{p}{\ga_0}
+\sum_{j=1}^qp^{2j-1}\ga_j\Big)\frac{h^{2\mu-2}}{p^{2s-2}}\norm{u}_{H^{s}(\Om)}^2.\nn
\end{align}
That is, \eqref{e4.6b} holds.

\eqref{e4.6d} can be proved similarly as above. It is clear
that \eqref{elem-app2-a} and \eqref{elem-app2-c} hold with
$s=2$ and \eqref{elem-app2-b} holds with $s=2$ and $j=1$,  that is,
\begin{align}
&\norm{u-\hat u_h}_{H^j(\Om)}
\lesssim \frac{h^{2-j}}{p^{2-j}}\norm{u}_{H^{2}(\Om)},
\quad j=0,1,\label{elem-app2-d}\db\\
&\norm{u-\hat u_h}_{H^2(\cT_h)}\lesssim p\norm{u}_{H^{2}(\Om)},
\label{elem-app2-e}\db\\
&\sum_{K\in\cT_h}\sum_{e\subset\pa K}\norml{\frac{\pa (u-\hat u_h)}{\pa n_e}}{e}^2
\lesssim h\norm{u}_{H^{2}(\Om)}^2.\label{elem-app2-f}
\end{align}
Since $u\in H^{q+1}(\Om)$, we have from Lemma~\ref{lem-inv}
and \eqref{elem-app2-e} that, for $2\le j\le q$,
\begin{align}\label{elem-app2-g}
\sum_{e\in\cE_h^I}&\norml{\jm{\frac{\pa^j (u-\hat u_h)}{\pa n_e^j}}}{e}^2
=\sum_{e\in\cE_h^I}\norml{\jm{\frac{\pa^j \hat u_h}{\pa n_e^j}}}{e}^2\\
&\lesssim\sum_{K\in\cT_h}\sum_{e\subset\pa K}\norml{\frac{\pa^j \hat u_h}{\pa n_e^j}}{e}^2
\lesssim\frac{p^2}{h}\norm{\hat u_h}_{H^j(\cT_h)}^2\nn\\
&\lesssim\Big(\frac{p^2}{h}\Big)^{2j-3}\norm{\hat u_h}_{H^2(\cT_h)}^2
\lesssim\Big(\frac{p}{h}\Big)^{2j-3}\,p^{2j-2}\norm{u}_{H^2(\cT_h)}^2.\nn
\end{align}
Now \eqref{e4.6d} follows by combining the equality in
\eqref{elem-app2-h} and \eqref{elem-app2-d}--\eqref{elem-app2-g}.
This completes the proof of the lemma.
\end{proof}

By combining Lemma~\ref{lem-4.2a} and Lemma~\ref{lem-app2}(i) we have the
following estimates for the projection error.
\begin{lemma}\label{lem-4.2}
Let $\mu= \min\set{p+1,s}$ and $q<\mu$.
Suppose problem \eqref{e1.1}--\eqref{e1.3} is $H^{s}$-regular
and \eqref{Areg} holds. Then there hold the following estimates:
\begin{align}
&\norm{u-\tuh^q}_{1,h,q} + \sqrt{\lambda k} \norml{u-\tuh^q}{\Gamma_R}
\lesssim \cerr{q}\,M_s(f,g)\,\frac{k^{s-1} h^{\mu-1}}{p^{s-1}}, \label{e4.2a}\\
&\norml{u-\tuh^q}{\Om}
\lesssim \frac{p}{k}\, \cerrh{q}\,M_s(f,g)\,\frac{k^s h^\mu}{p^s},  \label{e4.2b}
\end{align}
where
\begin{align*}
&\cerr{q}:= \lambda \Big(1+\frac{p}{\ga_0}+\sum_{j=1}^q p^{2j-1}\ga_j
+\frac{kh}{\lambda p}\Big)^{\frac12},\\
&\cerrh{q}:=\Big(1+\frac{1}{\ga_0\,p}+\frac{\gamma_1}{p}
+kh\Big)^{\frac12}\,\cerr{q},\qquad\lambda:=1+\frac{p}{\gamma_0}.
\end{align*}
\end{lemma}

\begin{remark}
The requirement $q<s$ in the lemma is clear since the projection
$\tuh^q$ is not defined for $q>s$. However, for $q<s$,
$\norm{u-\tuh^q}_{1,h,q}$ can be bounded without using full
regularity of $u$, and such a bound is also useful (see Lemma \ref{lem-5.1}).
\end{remark}

\subsection{Error estimates for $u-u_h^q$}
In this subsection we shall derive error estimates for scheme \eqref{edg}.
This will be done by exploiting the linearity of the Helmholtz equation
and making use of the stability estimates derived in Theorem \ref{thm_sta} and
the projection error estimates established in Lemma \ref{lem-4.2}.

Let $u$ and $u_h^q$ denote the solutions of  \eqref{e1.1}--\eqref{e1.3} and \eqref{edg},
respectively. Assume that $u\in H^s(\Ome)$ with $s\ge q+1$, then \eqref{2.6}  holds for
$v_h\in V_h^p$. Define the error function $e_h:=u-u_h^q$.
Subtracting \eqref{edg} from \eqref{2.6} yields the following error equation:
\begin{equation}\label{error-eq}
a_h^q(e_h,v_h) -k^2 (e_h,v_h) +\i k \langle e_h,v_h\rangle_{\Gamma_R} =0
\qquad \forall v_h\in V_h^p.
\end{equation}
Let $\tuh^q$ be the elliptic projection of $u$ as defined in the previous
subsection. Write $e_h=\eta-\xi$ with $\eta:=u-\tuh^q, \xi:=u_h^q-\tuh^q$.
From \eqref{error-eq} and \eqref{e4.4} we get
\begin{align}\label{e4.14}
a_h^q(\xi,v_h) -k^2 (\xi,v_h) +\i k \langle \xi,v_h\rangle_{\Gamma_R}
&= a_h^q(\eta,v_h) -k^2 (\eta,v_h) +\i k \langle \eta,v_h\rangle_{\Gamma_R} \\
&=-k^2 (\eta,v_h) \qquad \forall v_h\in V_h^p. \nn
\end{align}
The above equation implies that $\xi\in V_h^p$ is the solution of
scheme \eqref{edg} with sources terms $f=-k^2\eta$ and $g\equiv 0$.
Then an application of Theorem \ref{thm_sta} and Lemma \ref{lem-4.2} immediately gives

\begin{lemma}\label{lem-4.3}
$\xi=u_h^q-\tuh^q$ satisfies the following estimate:
\begin{align}\label{est_xi}
\norml{\xi}{\Om}+\frac{1}{k} \norm{\xi}_{1,h,q}
\lesssim \csta{q}\, k p \, \cerrh{q}\,M_s(f,g) \,\frac{k^s h^\mu}{p^s}.
\end{align}
\end{lemma}

We are ready to state our error estimate results for scheme \eqref{edg},
which follows from Lemma \ref{lem-4.3},  Lemma \ref{lem-4.2} and an
application of the triangle inequality.

\begin{theorem}\label{thm_main}
Let $u$ and $u_h^q$ denote the solutions of \eqref{e1.1}--\eqref{e1.3} and \eqref{edg},
respectively.  Let $\mu=\min\{p+1,s\}$ and $q<\mu$.
Assume assumption~\eqref{Areg} holds. Then
\begin{align}\label{e4.15}
&\norm{u-u_h^q}_{1,h,q} \lesssim
\Bigl( \cerr{q} + \csta{q} k^3 h\cerrh{q}\Bigr)\,M_s(f,g)\,\frac{k^{s-1} h^{\mu-1}}{p^{s-1}},\\
&\norml{u-u_h^q}{\Ome} \lesssim  p\,\cerrh{q}\, \Bigl( \frac{1}{k} + \csta{q} k \Bigr) \,M_s(f,g)\,\frac{k^s h^\mu}{p^s}.
\label{e4.16}
\end{align}
\end{theorem}

\begin{remark}
$q<s$ is required in the theorem because $\norm{u-u_h^q}_{1,h,q}$ is not defined
for $q>s$. However, we note that
the $hp$-IPDG solution $u_h^q$ is always well-defined regardless the
regularity of underlying PDE solution $u$. For $q<s$,
$\norm{u-u_h^q}_{1,h,q}$ can be bounded without using full
regularity of $u$, and such a bound is also useful (see Lemma \ref{lem-5.2}).
\end{remark}

By combining Theorem~\ref{thm_main} and Theorem~\ref{thm_sta_uniform}
we have the following theorem that gives the best convergence order
so far that we can obtain theoretically for the method \eqref{edg}
under the mesh condition $k^3h^2p^{-1}\gtrsim 1$ (cf. Theorem~\ref{thm_main2dg0}).

\begin{theorem}\label{thm_err0_uniform}
Under the assumptions of Theorem~\ref{thm_sta_uniform} and \ref{thm_main}, we have
\begin{align}\label{e4.17}
&\norm{u-u_h^q}_{1,h,q} +k\,\norml{u-u_h^q}{\Ome}\\
&\lesssim \left\{\begin{array}{ll}
k\,p^{\frac{8}{3}}h^{-\frac13}\,M_s(f,g)\,\frac{k^{s-1} h^{\mu-1}}{p^{s-1}}
&\text{ if } q=p,\\
k\,\max\set{p^{\frac{8}{3}}h^{-\frac13},p^{\frac{2q+4}{q+1}}h^{-\frac{1}{q+1}}}\,M_s(f,g)\,\frac{k^{s-1} h^{\mu-1}}{p^{s-1}} &\text{ if } q<p.
\end{array}
\right.\nn
\end{align}
\end{theorem}
\begin{proof}
The proof is obvious since $\cerr{q}\simeq\cerrh{q}\simeq 1.$
\end{proof}
\begin{remark}
(a) Estimates \eqref{e4.15}--\eqref{e4.17} are so-called preasymptotic error
estimates which are suboptimal in $h$ and $k$. They can be improved to 
optimal order when $k^3h^2p^{-1}\le C_0$, where $C_0$ is some constant 
independent of $k$, $h$, $p$, and the penalty parameters (see 
Theorem~\ref{thm_main2dg0}).  The second term on the right-hand side 
of \eqref{e4.15} is called a pollution
term for $\norm{u-u_h^q}_{1,h,q}$.

(b) Theorem~\ref{thm_err0_uniform} shows that
$\norm{u-u_h^q}_{1,h,q} +k\,\norml{u-u_h^q}{\Ome}\to 0$
if $q=p$ and $p^{\frac{11}{3}-s}k^s h^{\mu-\frac43}\to 0$,
or if $q=1<p$ and $p^{4-s}k^s h^{\mu-\frac32}\to 0$.
\end{remark}

%%%%%%%%%%%%%%%%%%
\section{Stability-error iterative improvement}\label{sec-SEII}
In this section we derive some improved optimal order stability and error estimates
for the $hp$-IPDG solution under the mesh condition that $k^3h^2p^{-1}\le C_0$  by using a stability-error iterative procedure, where $C_0$ is some constant independent of $k$, $h$, $p$, and the penalty parameters (see Theorem~\ref{thm_main2dg0}).

By combining  Lemma~\ref{lem-4.2a} and Lemma~\ref{lem-app2}(ii),
we have the following estimates for the projection error when only
the $H^2$-norm of the solution $u$ is allowed in the error bound.

\begin{lemma}\label{lem-5.1}
Suppose problem \eqref{e1.1}--\eqref{e1.3} is $H^{\max\set{q+1,2}}$-regular.
Then there hold the following estimates:
\begin{align}
&\norm{u-\tuh^q}_{1,h,q} + \sqrt{\lambda k} \norml{u-\tuh^q}{\Gamma_R}
\lesssim \cerr{2,q}\,M(f,g)\,\frac{k\,h}{p}, \label{e5.1a}\\
&\norml{u-\tuh^q}{\Om}
\lesssim \frac{p}{k}\, \cerrh{2,q}\,M(f,g)\,\frac{k^2 h^2}{p^2},  \label{e5.1b}
\end{align}
where
\begin{align*}
&\cerr{2,q}:= \lambda \Big(1+\frac{p}{\ga_0}+ p\,\ga_1+\sum_{j=2}^q p^{2j-2}\ga_j+\frac{kh}{\lambda p}\Big)^{\frac12},\\
&
\cerrh{2,q}:=\Big(1+\frac{1}{\ga_0\,p}+\frac{\gamma_1}{p}+kh\Big)^{\frac12}\,\cerr{2,q},\qquad\lambda:=1+\frac{p}{\gamma_0}.
\end{align*}
\end{lemma}
By a similar argument to that used to prove Theorem~\ref{thm_main},
we have the following error bounds which only involves $M(f,g)$.
\begin{lemma}\label{lem-5.2}
Let $u$ and $u_h^q$ denote the solutions of \eqref{e1.1}--\eqref{e1.3} and \eqref{edg},
respectively.  Suppose $u\in H^{\max\set{q+1,2}}(\Om)\cap H_{\Ga_D}^1(\Om)$. Then
\begin{align}\label{e5.2a}
&\norm{u-u_h^q}_{1,h,q} \lesssim
\Bigl( \cerr{2,q} + \csta{q} k^3 h\cerrh{2,q}\Bigr)\,M(f,g)\,\frac{k\,h}{p},\\
&\norml{u-u_h^q}{\Ome} \lesssim  p\,\cerrh{2,q}\, \Bigl( \frac{1}{k} + \csta{q} k \Bigr) \,M(f,g)\,\frac{k^2 h^2}{p^2}.
\label{e5.2b}
\end{align}
\end{lemma}

We are now ready to state our final main theorem of this paper.

\begin{theorem}\label{thm_main2dg0}
Let $u$ and $u_h^q$ denote the solutions of \eqref{e1.1}--\eqref{e1.3} and \eqref{edg},
respectively.  Suppose $u\in H^{\max\set{q+1,2}}(\Om)\cap H_{\Ga_D}^1(\Om)$.  Assume that $k\gtrsim 1$, $k\,h\lesssim 1$, and that
$p{\ga_0}^{-1}+\sum_{j=1}^q p^{2j-1}\ga_j\lesssim 1$. Then there exists
a constant $C_0>0$, which is independent of $k$, $h$, $p$, and the penalty parameters,
such that if $k^3 h^2p^{-1}\le  C_0$, then the following stability
estimates hold:
\begin{align}\label{esta1}
&\norm{u_h^q}_{1,h,q}\lesssim  M(f,g),\\
&\norml{u_h^q}{\Ome}\lesssim  \frac{1}{k} M(f,g). \label{esta1a}
\end{align}
Moreover, under the assumption~\eqref{Areg}, there hold the following error estimates:
\begin{align}\label{eerr1}
&\norm{u-u_h^q}_{1,h,q} \lesssim ( 1 + k^2 h) \,M_s(f,g)\,\frac{k^{s-1} h^{\mu-1}}{p^{s-1}},\\
&\norml{u-u_h^q}{\Ome} \lesssim \,p\,M_s(f,g)\,\frac{k^s h^\mu}{p^s}. \label{eerr1a}
\end{align}
\end{theorem}
\begin{proof} We only prove \eqref{esta1} since \eqref{esta1a} can be proved
similarly and \eqref{eerr1}--\eqref{eerr1a} follow from the improved stability
estimates and the argument used in the proof of Theorem~\ref{thm_main}.

From Theorem~\ref{thm_sta} we have
\begin{equation}\label{e5.1}
\norm{u_h^q}_{1,h,q}\lesssim k\,\csta{q}\, M(f,g),
\end{equation}
where $\csta{q}$ is defined in \eqref{csta}. From Lemma~\ref{lem-5.2} we have
\begin{equation}\label{e5.2}
    \norm{u-u_h^q}_{1,h,q} \lesssim
\Bigl( \cerr{2,q} + k\,\csta{q} k^2 h\cerrh{2,q}\Bigr) \frac{k\,h}{p}\,M(f,g).
\end{equation}
Now it follows from Theorem~\ref{stability} and the triangle inequality that
\begin{align}\label{e5.3}
\norm{u_h^q}_{1,h,q} &\le\norm{u}_{1,h,q}+\norm{u-u_h^q}_{1,h,q}
=|u|_{1,h}+\norm{u-u_h^q}_{1,h,q}\\
&\lesssim \Bigl(1+\cerr{2,q}\,\frac{k\,h}{p}
+ \frac{k^3 h^2}{p}\cerrh{2,q}\,k\,\csta{q}  \Bigr) M(f,g).\nn
\end{align}
Repeating the above process yields that there exists a constant $C_1$
independent of $k$, $h$, $p$, and the penalty parameters, and a sequence of
positive numbers $\La_j$ such that
\begin{align}\label{e5.4}
\norm{u_h^q}_{1,h,q}\le \La_jM(f,g),
\end{align}
with
\begin{align*}
\La_0\simeq k\,\csta{q},\quad \La_j
=C_1(1+\cerr{2,q}\,\frac{k\,h}{p}) +  C_1\,\cerrh{2,q}\,\frac{k^3 h^2}{p}\,\La_{j-1},
\quad j=1,2,\cdots.
\end{align*}
A simple calculation yields that if  $C_1\,\cerrh{2,q}\,k^3 h^2<\ta\,p$ for
some positive constant $\ta<1$ then
\[
\lim_{j\to\infty}\La_j=\frac{C_1(p+\cerr{2,q}\,k\,h) }{p-C_1\,\cerrh{2,q}\,k^3 h^2},
\]
which implies \eqref{esta1} by noting that $\cerr{2,q}, \cerrh{2,q}\lesssim 1$
and that $\cerr{2,q}\,k\,h\lesssim \big(\cerr{2,q}^2k^3 h^2\big)^{\frac12}
\lesssim \big(\cerrh{2,q}^2k^3 h^2\big)^{\frac12}\lesssim p^{\frac12}$.
\end{proof}

Note that the stability estimates in \eqref{esta1} and \eqref{esta1a} are of the same order as the PDE stability estimates given
in Theorem~\ref{stability}.

%%%%%%%%%%%%%%

\end{document}